\documentclass[11pt,a4paper]{elsarticle}

\makeatletter
\def\ps@pprintTitle{%
\let\@oddhead\@empty
\let\@evenhead\@empty
\def\@oddfoot{}%
\let\@evenfoot\@oddfoot}
\makeatother

\usepackage[utf8]{inputenc}

\usepackage{amsmath}
\usepackage{amsthm}
\usepackage{enumerate}
\usepackage{mathdots}
\usepackage{amsfonts}
\usepackage{multicol}
\usepackage{amsmath}
\usepackage{amssymb}
\usepackage{fancybox}
\usepackage[usenames,dvipsnames]{color}
\usepackage{graphicx}
\usepackage{tikz}
\usepackage[left=3cm,top=4cm,right=3cm,bottom=3.5cm]{geometry}
\usepackage{arydshln}
\usepackage{ulem}

\DeclareMathOperator{\orb}{\mathcal{O}}

\newcommand{\NN}{\mathbb{N}}

\newcommand{\CC}{\mathbb{C}}

\newcommand{\FF}{\mathbb{F}}

\newcommand{\la}{\lambda}

\def\rank{\mathop{\rm rank}\nolimits}

\newcommand{\blue}{\color{black}}

\newtheorem{theo}{Theorem}[section]

\newtheorem{deff}[theo]{Definition}

\newtheorem{lem}[theo]{Lemma}

\newtheorem{cor}[theo]{Corollary}

\newtheorem{rem}[theo]{Remark}

\newtheorem{example}[theo]{Example}

\begin{document}

	\title{Minimal rank factorizations of polynomial matrices}

	\author[orebro]{Andrii Dmytryshyn\fnref{fn2}}
	\ead{andrii.dmytryshyn@oru.se}

    \author[uc3m]{Froil\'{a}n Dopico\fnref{fn1}}
	\ead{dopico@math.uc3m.es}

	\author[ucl]{Paul~Van~Dooren}
	\ead{paul.vandooren@uclouvain.be}
	
\address[orebro]{Department of Mathematical Sciences, Chalmers University of Technology and University of Gothenburg, 41296 Gothenburg, Sweden.}
	
		\fntext[fn2]{Supported by the Vetenskapsr\aa det (Swedish Research Council)
[grant 2021-05393].}
	
	\address[uc3m]{Departamento de Matem\'aticas,
		Universidad Carlos III de Madrid, Avda. Universidad 30, 28911 Legan\'es, Spain.}

	\fntext[fn1]{Supported {\blue by the ``Agencia Estatal de Investigaci\'on'' of Spain MCIN/AEI/ 10.13039/501100011033} through grants PID2019-106362GB-I00, {\blue PID2023-147366NB-I00 and  RED2022-134176-T,} and by the Madrid Government (Comunidad de Madrid-Spain) under the Multiannual Agreement with UC3M in the line of Excellence of University Professors (EPUC3M23), and in the context of the V PRICIT (Regional Programme of Research and Technological Innovation).}

		\address[ucl]{Department of Mathematical Engineering, Universit\'{e} catholique de Louvain, Avenue Georges Lema\^itre 4, B-1348 Louvain-la-Neuve, Belgium.}
		

\begin{abstract}
We investigate rank revealing factorizations {\blue of $m \times n$} polynomial matrices $P(\la)$ into products of three, $P(\la) = L(\la) E(\la) R(\la)$, or two, $P(\la) = L(\la) R(\la)$, polynomial matrices. Among all possible factorizations of these types, we focus on those for which $L(\la)$ and/or $R(\la)$ is a minimal basis, since they {\blue have favorable properties from the point of view of data compression and} allow us to relate easily the degree of $P(\la)$ with some degree properties of the factors. We call these factorizations {\it minimal rank factorizations}. Motivated by the well-known fact that, generically, rank deficient polynomial matrices over the complex field do not have eigenvalues, we pay particular attention to the properties of the minimal rank factorizations of polynomial matrices without eigenvalues. We carefully analyze the degree properties of generic minimal rank factorizations in the set of complex $m \times n$ polynomial matrices with normal rank at most {\blue $r< \min \{m,n\}$} and degree at most $d$, and we prove that {\blue there are only $rd+1$ different classes of generic factorizations according to the degree properties of the factors and that all of them} are of the form $L(\la) R(\la)$, where the degrees of the $r$ columns of $L(\la)$ differ at most by one, the degrees of the $r$ rows of $R(\la)$ differ at most by one, and, for each $i=1, \ldots, r$, the sum of the degrees of the $i$th column of $L(\la)$ and of the $i$th row of $R(\la)$ is equal to $d$. Finally, we show how these sets of polynomial matrices with generic factorizations are related to the sets of polynomial matrices with generic eigenstructures.
\end{abstract}

	\begin{keyword} polynomial matrix \sep factorization \sep
		complete eigenstructure\sep genericity\sep minimal bases \sep normal rank
		
		\medskip\textit{AMS subject classifications}: 15A18, 15A22, 15A23, 15A54
	\end{keyword}

	\maketitle
	
\section{Introduction} \label{sec.intro}
Given an $m\times n$ matrix $A$ with complex entries and {\blue rank $r$}, it is often useful to express $A$ as the product of three factors $A = L E R$ of sizes $m\times r$, $r\times r$ and $r\times n$, respectively, or as the product of two factors $A = L R$ of sizes $m\times r$ and $r\times n$, respectively. Such factorizations are sometimes called rank-revealing factorizations, or rank factorizations for short,  since the sizes of the factors reveal the rank of the matrix. The singular value decomposition is probably the best known example of a rank-revealing factorization, though several other rank-revealing factorizations exist and are used in practice. Rank-revealing factorizations have many applications. Among them, data compression {\blue when $r \ll \min \{m , n\}$} plays an important role \cite{Golub-VanLoan}. {\blue Another relevant application is to use  rank-revealing factorizations as compact representations, or parametrizations, of the elements in the manifold of $m\times n$ matrices with rank $r$ that allow for the efficient solution of some optimization problems on this manifold \cite[Secs. 2.6, 2.8, 7.5]{boumalbook}.} It is well-known that a rank-revealing factorization $A = L R$ is equivalent to expressing $A$ as a sum of $r$ rank-1 matrices $A = v_1 u_1^T + \cdots + v_r u_r^T$, where $v_1, \ldots, v_r$ are the columns of $L$ and $u_1^T , \ldots , u_r^T$ are the rows of $R$. {\blue The concepts mentioned in this introduction are revised in Section \ref{sec.prelim}.}

The main goal of this paper is to investigate rank-revealing factorizations of $m\times n$ polynomial matrices $P(\la)$, of normal rank {\blue $r$} and degree $d$, into products of three, $P(\la) = L(\la) E(\la) R(\la)$, or two, $P(\la) = L(\la) R(\la)$, {\it polynomial matrices}. We will see that this problem is very different from the corresponding one for constant matrices and that it requires the use of completely different tools. These differences come essentially from two facts. First, from the constraint that the factors must be also polynomial matrices and, second, from the notion of degree, and the non-trivial question of how the degree of $P(\la)$ is related to the degrees (of the entries) of the factors. {\blue For instance, the naive idea that the sum of the degrees of the polynomial factors is equal to the degree of $P(\la)$ is not valid because there exist polynomial matrices for which none of their rank-reveling factorizations satisfy such relation (see, for instance, Example \ref{ex.arbitrarydeg} and Lemma \ref{lemm:lowbounddeg}).  Observe in this respect that the degrees of the factors play a key role for data compression since the matrix coefficients of $P(\la) = \sum_{i=0}^d P_i \la^i$ require to store $(d+1) m n$ numbers, while to store the coefficients of the factors in a rank-revealing factorization $P(\la) = L(\la) R(\la)$ requires to store up to $(d_L +1) mr + (d_R +1)nr$ numbers, where $d_L$ and $d_R$ are the degrees of $L(\la)$ and $R(\la)$ respectively. Clearly, high values of $d_L$ and $d_R$ are not desirable in terms of data compression. These difficulties extend to the possible use of rank-revealing factorizations of polynomial matrices as compact representations or parametrizations of the elements of the set of polynomial matrices of degree $d$ and normal rank $r$, because it is not clear which are the possible degrees that can be assigned to each factor and, so, it is not clear how to develop economic representations or how many of them are needed.}

{\blue These degree problems motivate us to focus on rank-revealing factorizations of polynomial matrices where $L(\la)$ is a minimal basis \cite{For75} of the column space of $P(\la)$ and/or $R(\la)$ is a minimal basis of the row space of $P(\la)$. We call these factorizations {\it minimal rank factorizations}. We prove that these factorizations have three important advantages. In the first place, the very same definition of minimal bases in \cite{For75} imply that they are the most economical bases of a rational subspace in terms of data storage. Combining this property with the fact that for any rank-revealing factorization of $P(\la)$, $L(\la)$ is a basis of the column space of $P(\la)$ and $R(\la)$ is a basis of the row space of $P(\la)$ (see Lemma \ref{lemm:fact0}), we see that minimal  rank factorizations are optimal in terms of data compression. This property is particularly relevant in the important generic case of rank deficient polynomial matrices without eigenvalues (see Theorem \ref{thm.minfact}), since in this case a middle factor $E(\la)$ is not necessary. The second key property of minimal rank factorizations is that they allow us to relate in a very clear way the degree of $P(\la)$ with certain matching properties of the degrees of the entries of the factors (see Corollary \ref{thm.degreepredic} in general and Theorem \ref{thm.minfact} for polynomial matrices without eigenvalues). Finally, in the case of polynomial matrices $P(\la)$ with eigenvalues, their minimal rank factorizations with three factors and with $L(\la)$ and $R(\la)$ both minimal bases guarantee that the middle factor $E(\la)$ contains all the finite eigenvalues of $P(\la)$ with their partial multiplicities (see Theorem \ref{thm.minfact0}-(i) and Remark \ref{rem.uniqueminfact}). Observe that if $r \ll \max \{ m, n\}$, the $r\times r$ factor $E(\la)$ is much smaller than $P(\la)$.

Despite the advantages described in the previous paragraph, minimal rank factorizations still allow for a lot of freedom on the possible degrees of the columns of $L(\la)$ and the rows of $R(\la)$, except when $d$ and $r$ are both small (see the discussion at the end of Section \ref{sec.minfactorizations}). Thus, further work is needed for finding a small number of compact parametrizations via factorizations of (a dense subset of) the set of polynomial matrices of degree $d$ and normal rank $r$ that might potentially allow us, for instance, to solve efficiently optimization problems on this set. Such compact parametrizations are related to some results available in the literature which are discussed in the next paragraph.}

{\blue It is well-known \cite{dmy-dop-2017} that {\it generic} $m\times n$ polynomial matrices with normal rank $r < \min \{m , n \}$ and degree at most $d$, over the complex field, do not have eigenvalues and have minimal indices with very particular properties. More precisely, the $m-r$ left minimal indices differ at most by one and the same happens with the $n-r$ right minimal indices. These properties combined with the Index Sum Theorem \cite{DDM} (see also Theorem \ref{thm.indexsum} below) give rise to the existence of only $rd+1$ different generic complete eigenstructures in the set $\CC[\la]^{m\times n}_{d,r}$ of $m\times n$ complex polynomial matrices with normal rank at most $r$ and degree at most $d$, and allow us to express $\CC[\la]^{m\times n}_{d,r}$ as the union of the closures of the $rd+1$ sets (usually called orbits) of all polynomial matrices with such generic eigenstructures (see Theorem \ref{thm.mainthandrii2017} below).

The generic results described above motivate us, in the first place, to study in detail the minimal rank factorizations of polynomial matrices without eigenvalues and, in the second place, to look for alternative descriptions of the set $\CC[\la]^{m\times n}_{d,r}$ in terms of the union of the closures of a few sets of polynomial matrices which have generic rank-revealing factorizations with very specific properties, instead of in terms of a few generic eigenstructures. In this line, we prove that, generically, {\it the polynomials in $\CC[\la]^{m\times n}_{d,r}$ can be factorized in only $rd+1$ different ways} according to the degree properties of the factors. More precisely, among other results, we prove that generically for $P(\la) \in  \CC[\la]^{m\times n}_{d,r}$ a factorization of $P(\la) = L(\la) R(\la)$ into two polynomial matrices of sizes $m\times r$ and $r\times n$ satisfies that the degrees of the $r$ columns of $L(\la)$ differ at most by one, the degrees of the $r$ rows of $R(\la)$ differ at most by one, and, for each $i = 1, \ldots , r$, the sum of the degrees of the $i$th column of $L(\la)$ and of the $i$th row of $R(\la)$ is equal to $d$, and that there are only $rd+1$ different ways to choose the involved column and row degrees (see Definition \ref{def.homogsubsets} and Theorem \ref{thm.4descrition}-(iii) among other results in this spirit). We emphasize that each of such $rd+1$ sets of polynomial matrices with these specific factorizations can be easily and efficiently parameterized using the vector coefficients of the columns of $L(\la)$ and the rows of $R(\la)$. In this context, we also study how the orbits of the polynomial matrices with the generic eigenstructures identified in \cite{dmy-dop-2017} (see also Theorem \ref{thm.mainthandrii2017}) are related to the polynomial matrices with the generic factorizations that we identify in this work (see Theorem \ref{thm.relationgeneigen}). These two generic views of the set $\CC[\la]^{m\times n}_{d,r}$ are complementary. The generic orbits provide geometrical insights into the structure of $\CC[\la]^{m\times n}_{d,r}$. But in order to use such insights for solving certain problems numerically, e.g., optimization or distance problems related to low rank polynomial matrices with given degree, parametrizations of these orbits are needed. The generic factorizations presented in this work provide such a parametrizations for the closures of the generic orbits. In summary, these factorizations help us to bridge the geometry and numerics for low rank matrix polynomials of given degree.
}

We are not aware of other similar results available in the literature, dealing with rank-revealing factorizations of polynomial matrices of degree larger than one. However, there exist factorizations of this type in the case of degree at most one, that is, in the case of matrix pencils. In fact, rank-revealing factorizations expressed as the sum of matrix pencils with rank one exist for unstructured pencils \cite{DDSIMAX2007,DDSIMAX2016,DDL} and also for matrix pencils with symmetry structures \cite{DMMFound}. {\blue We will explore in Remark \ref{rem.compen1}, Example \ref{ex.compen1}, and Remark \ref{rem.compen2} the relationship between the results previously obtained for unstructured pencils and the new ones for unstructured polynomial matrices of degree larger than one developed in this work.

It is worth to point out that} rank-revealing factorizations of matrix pencils have played a fundamental role in the study of the generic effect of low rank perturbations on the eigenstructure of a given regular matrix pencil. {\blue More precisely, in \cite{DDSIMAX2016} the authors express the set of all possible perturbations, i.e., $\CC[\la]^{n\times n}_{1,r}$, as the union of $r+1$ sets of pencils whose elements have rank-revealing factorizations with very specific properties (see \cite[Lemma 3.1]{DDSIMAX2016} or \cite[Lemma 4]{DDL}). After that, they define $r+1$ surjective smooth maps $\Phi_0 , \Phi_1, \ldots, \Phi_r$ from $\CC^{3rn}$ onto each of these sets using the specific forms of the factors of the elements of these sets \cite[Definition 3.2]{DDSIMAX2016}. Finally, it is proved in \cite[Theorem 3.4]{DDSIMAX2016} that there exist $r+1$ generic subsets $G_0, G_1, \ldots , G_r$ of $\CC^{3rn}$ such that all the perturbation pencils in the sets $\Phi_0(G_0), \Phi_1(G_1), \ldots , \Phi_r(G_r)$ produce the same ``generic'' effect on the eigenstucture of the unperturbed regular pencil. We hope that this strategy combined with the results developed in this paper about expressing $\CC[\la]^{n\times n}_{d,r}$ as the union of the closures of $rd+1$ sets of polynomial matrices having rank-revealing factorizations with the very specific properties described above will have applications in the study of the generic effect of low rank perturbations on the eigenstructure of a given regular {\it polynomial matrix} of degree larger than one, which is a problem that remains open in the literature.}

We emphasize that the rank-revealing factorizations of matrix pencils {\blue in \cite{DDSIMAX2007,DDSIMAX2016,DDL,DMMFound}} have been obtained by using the Kronecker canonical form of pencils under strict equivalence \cite{Gan59}, or structured versions of this form. Since a canonical form of this type does not exist for polynomial matrices of degree larger than one, the problem for polynomial matrices is harder than for matrix pencils, {\blue requires different tools, and yields results weaker than those for pencils (see the discussions in Remark \ref{rem.compen1}, Example \ref{ex.compen1}, and Remark \ref{rem.compen2})}.

The paper is organized as follows. Section \ref{sec.prelim} includes some known concepts and results that are important for obtaining the main results of this paper. Rank-revealing factorizations and minimal rank factorizations of polynomial matrices are introduced in Section \ref{sec.minfactorizations}, where their properties are also studied. Section \ref{sect.compact} establishes the generic properties of rank-revealing factorizations and minimal rank factorizations. Finally, Section \ref{sec.conclusions} presents some conclusions and possible lines of future research.

\section{Preliminaries} \label{sec.prelim} This section summarizes the notation and some of the results previously published in the literature, that will be used in the paper. Many of the results in this paper are valid over an arbitrary field $\FF$ while others are only valid over the field $\CC$ of complex numbers. This will be clearly indicated in the text by using either $\FF$ or $\CC$. $\FF[\la]$ stands for the ring of polynomials in the variable $\la$ with coefficients in $\FF$ and $\FF(\la)$ stands for the field of fractions of $\FF[\la]$, i.e., rational functions in the variable $\la$ with coefficients in $\FF$. A polynomial vector is a vector with entries in $\FF[\la]$. $\FF[\la]^{m\times n}$ and $\FF(\la)^{m\times n}$ denote the sets of $m\times n$ polynomial matrices and of $m\times n$ rational matrices, respectively, over $\FF$. The degree of a polynomial vector, $q(\la)$, or of a polynomial matrix, $P(\la)$, is the highest degree of all of its entries and is denoted by $\deg(q)$ or $\deg(P)$. The degree of the zero polynomial is defined to be $-\infty$. The set of $m\times n$ polynomial matrices of degree {\it at most $d$} is denoted by $\FF[\la]^{m\times n}_d$. Given a list $\underline{\mathbf{d}} = (d_1, d_2, \ldots ,d_m)$ of nonnegative integers, $\FF[\la]^{m\times n}_{\underline{\mathbf{d}}}$ denotes the set of $m\times n$ polynomial matrices whose $i$th row has degree {\it at most} $d_i$ for $i =1, \ldots , m$.
We also use $\overline \FF$ for the algebraic closure of $\FF$, $I_n$ for the $n\times n$ identity matrix, and $0_{m\times n}$ for the $m\times n$ zero matrix, where the sizes are omitted when they are clear from the context. We need to use very often the $i$th row or the $j$th column of a polynomial matrix $P(\la)$ and we adopt the following compact notations for them: $P_{i*} (\la)$, or simply $P_{i*}$, denotes the $i$th row of $P(\la)$ and $P_{*j} (\la)$, or simply $P_{*j}$, denotes the $j$th column of $P(\la)$.

The normal rank of a polynomial or rational matrix $P(\la)$, denoted as $\mbox{rank} (P)$, is the rank of $P(\la)$ considered as a matrix over the field $\FF (\la)$, or the size of the largest non-identically zero minor of $P(\la)$. The reader can find more information on polynomial and rational matrices in the books \cite{Gan59,Kai80}.

The set of $m \times n$ polynomial matrices with degree {\it at most} $d$ and normal rank {\it at most} $r$ is denoted by $\FF[\la]^{m\times n}_{d,r}$. In the case $\FF = \CC$ and $r < \min \{m,n\}$, new results about factorizations of the elements of this set will be presented in Section \ref{sect.compact}. In order to avoid trivialities, every time that the symbol $\FF[\la]^{m\times n}_{d,r}$ is written it should be understood that the integers $d$ and $r$ satisfy $d \geq 1$ and $r\geq 1$.

The well-known Smith form of a polynomial matrix plays a very important role in this work and the corresponding result is presented in Theorem \ref{thm.smith} \cite{Gan59}. It requires the use of unimodular polynomial matrices, that is, square polynomial matrices with constant nonzero determinant.

\begin{theo} {\rm (Smith form)} \label{thm.smith}
Let $P(\la) \in \FF[\la]^{m\times n}$ with $\rank (P) = r$. Then there exist {\blue a unique diagonal matrix $S(\la) \in \FF[\la]^{m\times n}$ and  unimodular matrices $U(\la) \in \FF[\la]^{m\times m}$, $V(\la) \in \FF[\la]^{n\times n}$} such that
\begin{equation} \label{eq.snf}
    P(\lambda)=   U(\lambda)S(\lambda)V(\lambda), \quad
    S(\la):= \left[
      \begin{array}{c|c}
        \begin{array}{cccc}
        e_1(\lambda) & 0 & \ldots & 0 \\
        0 & e_2(\lambda) & \ddots & \vdots \\
        \vdots & \ddots & \ddots & 0  \\
        0 & \ldots & 0 & e_r(\lambda)
        \end{array}
        &  0_{r \times (n-r)} \\[-8pt] & \\ \hline
        0_{(m-r) \times r} & 0_{(m-r) \times (n-r)}
\end{array} \right],
\end{equation}
  where each polynomial $e_j(\lambda) \in \FF [\la]$ is monic and divides $e_{j+1}(\lambda)$ {\blue for $j=1,\ldots,r-1$.}
\end{theo}

The unique matrix $S(\la)$ in \eqref{eq.snf} is the Smith form of $P(\la)$ and the expression $P(\lambda)=   U(\lambda)S(\lambda)V(\lambda)$ is called a Smith factorization of $P(\la)$. Smith factorizations are not unique. The polynomials $e_j(\lambda)$ are called the invariant polynomials of $P(\lambda)$ and those that are equal to $1$ are called trivial invariant polynomials. For any $\alpha \in \overline{\FF}$, the invariant polynomials can be uniquely factorized as
$e_j (\la) = (\la - \alpha)^{\sigma_j} p_j(\la)$, with $p_j(\la) \in \overline{\FF} [\la]$, $p_j(\alpha) \ne 0$ and $\sigma_j \in \NN =\{0, 1, 2, \ldots \}$, for $j=1, \ldots ,r$. The sequence $\sigma_1 \leq \cdots \leq \sigma_r$ is called the partial multiplicity sequence of $P(\la)$ at $\alpha$. A root $\beta \in \overline{\FF}$ of any of the invariant polynomials $e_j (\la)$ of $P(\la)$ is called a finite eigenvalue of $P(\la)$. Equivalently, $\beta \in \overline{\FF}$ is a finite eigenvalue of $P(\la)$ if and only if the partial multiplicity sequence of $P(\la)$ at $\beta$ contains at least one nonzero term.

The partial multiplicity sequence at $\infty$ of $P(\la) \in \FF[\la]_d^{m\times n}$ is defined to be the partial multiplicity sequence at $0$ of $\la^d P(1/\la) \in \FF[\la]_d^{m\times n}$ and it is said that $P(\la)$ has an eigenvalue at $\infty$ if its partial multiplicity sequence at $\infty$ contains at least one nonzero term, or, equivalently, if zero is an eigenvalue of $\la^d P(1/\la)$. It is easy to prove that the first term of the partial multiplicity sequence at $\infty$ and the degree of {\blue the polynomial matrix} are related as follows.
\begin{lem} \label{lemm.degandinfty} {\rm \cite[Lemma 2.6]{quasi-trig}} Let $P(\la) \in \FF[\la]_d^{m\times n}$ with $\rank (P) =r$ and partial multiplicity sequence at $\infty$ equal to $0 \leq \gamma_1 \leq \gamma_2 \leq \cdots \leq \gamma_r$. Then $\gamma_1 = d - \deg(P)$.
\end{lem}

{\blue We remark that if a polynomial matrix $P(\la) \in \FF[\la]_d^{m\times n}$, then $P(\la) \in \FF[\la]_e^{m\times n}$ for any $e>d$. Thus, the definition above of the partial multiplicity sequence at $\infty$ and of eigenvalue at $\infty$ of $P(\la)$ depends on the choice of the set to which $P(\la)$ belongs, though such dependence is trivial via the shift $e-d$ of the sequence. Therefore, every time we mention the eigenvalue at $\infty$ of a polynomial matrix, we will specify the set $\FF[\la]_d^{m\times n}$ to which the polynomial belongs. This dependence of the partial multiplicity sequence at $\infty$ on the chosen set $\FF[\la]_d^{m\times n}$ is rather extensively discussed in the literature. In fact, the set $\FF[\la]_d^{m\times n}$ is what is called in \cite{mobiusmackeys} the vector space of $m\times n$ polynomial matrices of {\it grade} $d$, its elements are said to have grade $d$, independently of the degree they may have, and their partial multiplicity sequences at $\infty$ are defined with respect to $d$. The grade and the corresponding definition of eigenvalue at $\infty$ have been used very often for solving several problems on polynomial matrices as, for instance, genericity problems \cite{DeDD24,dmy-dop-2017,dmy-dop-2018}, the analysis of certain structured polynomial matrices \cite{smithpalin}, and in M\"{o}bius and more general rational transformations \cite{mobiusmackeys,noferinirational}.}

Next, we {\blue recall} the concept of minimal bases of a rational subspace \cite{For75}. Let us consider the vector space $\FF(\la)^n$ over the field $\FF(\la)$. A subspace $\mathcal{V}$ of $\FF(\la)^n$ is called a rational subspace. {\blue By clearing out denominators, one can} see that every rational subspace $\mathcal{V}$ has bases consisting entirely of polynomial vectors, {\blue which are called polynomial bases of $\mathcal{V}$.} Following Forney \cite{For75}, we say that a minimal basis of $\mathcal{V}$ is a {\blue polynomial} basis of $\mathcal{V}$ consisting of vectors whose sum of degrees is minimal among all {\blue polynomial} bases of $\mathcal{V}$. A key property \cite{For75} is that the ordered list of degrees of the polynomial vectors in any minimal basis of $\mathcal{V}$ is always the same. These degrees are called the minimal indices of $\mathcal{V}$. {\blue Observe that if $\dim \mathcal{V} = p$ and the degrees of the vectors in a polynomial basis $\mathcal{B}$ of $\mathcal{V}$ are $d_1, d_2, \ldots, d_p$, then $n \sum_{i=1}^p (d_i + 1)$ scalars of $\FF$ are needed to store $\mathcal{B}$. Therefore, the minimal bases of $\mathcal{V}$ are optimal from the point of view of data storage among the polynomial bases of $\mathcal{V}$, as we pointed out in Section \ref{sec.intro}. Minimal bases and indices also satisfy the following  ``Strong Minimality Property of Minimal Indices''.

\begin{theo} \label{thm.mackeystrong} {\rm \cite[Thm. 4.2]{Mackey-fil}} Let $\mathcal{V} \subseteq \FF(\la)^n$ be a $p$-dimensional rational subspace with minimal indices $\varepsilon_1 \leq \varepsilon_2 \leq \cdots \leq \varepsilon_p$ and let  $d_1 \leq d_2 \leq \cdots \leq d_p$ be the ordered degrees of the vectors in a polynomial basis $\mathcal{B}$ of $\mathcal{V}$. Then, $\varepsilon_i \leq d_i$ for $i=1, \ldots, p$.
\end{theo}
}

The minimal bases of any rational subspace can be characterized in different important ways \cite[p. 495]{For75} (see also \cite{Kai80}). Among them, we emphasize the characterization in Theorem \ref{thm.minbasischar}, which {\blue leads us to introduce} Definition \ref{def.colred}.

\begin{deff}\label{def.colred} {\rm {\blue \cite[Def. 2.5.6, p. 27]{wolovich}}}
Let $d'_1,\ldots,d'_n$ be the degrees of the columns of $N(\lambda) \in \FF[\la]^{m\times n}$.
The highest-column-degree coefficient matrix of $N(\la)$, denoted by $N_{hc}$,
is the $m\times n$ constant matrix whose {\rm $j$th} column
is the vector coefficient of $\lambda^{d'_j}$ in the {\rm $j$th} column of $N(\lambda)$. The polynomial matrix $N(\la)$ is said to be column reduced if $N_{hc}$ has full column rank.

Similarly, let $d_1,\ldots,d_m$ be the degrees of the rows of $M(\la) \in \FF[\la]^{m\times n}$. The highest-row-degree coefficient matrix of $M(\la)$, denoted by $M_{hr}$, is the $m\times n$ constant matrix whose {\rm $j$th} row is the vector coefficient of $\lambda^{d_j}$ in the {\rm $j$th} row of $M(\lambda)$.
The polynomial matrix $M(\la)$ is said to be row reduced if $M_{hr}$ has full row rank.
\end{deff}

\begin{theo}\label{thm.minbasischar} {\rm {\blue \cite[Main Thm. 2, p. 495]{For75}}}
  The columns {\rm (}resp., rows{\rm )} of a polynomial matrix $N(\lambda) \in \FF[\la]^{m\times n}$ are a minimal basis of the rational subspace they span
  if and only if $N(\lambda_0)$ has full column {\rm (}resp., row{\rm )} rank
  for all $\lambda_0 \in \overline \FF$,
  and $N(\la)$ is column {\rm (}resp., row{\rm )} reduced.
\end{theo}

Next, we define four rational subspaces associated with a polynomial matrix $P(\la)$.

\begin{deff} \label{def.ratsubsP} {\rm (Rational subspaces of a polynomial matrix)} Let $P(\la) \in \FF[\la]^{m\times n}$. Then
\begin{itemize}
  \item[\rm (i)] $\mathcal{N}_\ell (P) = \{y(\la) \in \FF(\la)^{1\times m}  : \, y(\la) P(\la) =0\} \subseteq \FF(\la)^{1\times m}$ is the left nullspace of $P(\la)$,
  \item[\rm (ii)] $\mathcal{N}_r (P) = \{x(\la) \in \FF(\la)^{n\times 1} : \, P(\la) x(\la) =0\} \subseteq \FF(\la)^{n\times 1}$ is the right nullspace of $P(\la)$,
  \item[\rm (iii)] ${\mathcal Row} (P) = \{w(\la) P(\la) \, : \, w(\la) \in \FF(\la)^{1\times m}\} \subseteq \FF(\la)^{1\times n}$ is the row space of $P(\la)$,
  \item[\rm (iv)] ${\mathcal Col} (P) = \{P(\la) v(\la) \, : \, v(\la) \in \FF(\la)^{n\times 1}\} \subseteq \FF(\la)^{m\times 1}$ is the column space of $P(\la)$.
\end{itemize}
\end{deff}

Observe that if $\rank (P) = r$, then $\dim \mathcal{N}_\ell (P) = m-r$, $\dim \mathcal{N}_r (P) = n-r$ and $\dim {\mathcal Row} (P) = \dim {\mathcal Col} (P) = r$, by the rank-nullity theorem \cite[Vol. I, p. 64]{Gan59}. Thus, $\mathcal{N}_\ell (P)$ has $m-r$ minimal indices, $\mathcal{N}_r (P)$ has $n-r$ minimal indices, and ${\mathcal Row} (P)$ and ${\mathcal Col} (P)$ have each of them $r$ minimal indices.

Given a polynomial matrix $P(\la) \in \FF[\la]^{m\times n}_d$, the set formed by its invariant polynomials, by its partial multiplicity sequence at $\infty$, by the minimal indices of $\mathcal{N}_\ell (P)$ and by the minimal indices of $\mathcal{N}_r (P)$ is often called the complete eigenstructure of $P(\la)$ \cite{DDVD,VDoorenDewilde}. Observe that the minimal indices of ${\mathcal Row} (P)$ and ${\mathcal Col} (P)$ are {\it not} included in the complete eigenstructure of $P(\la)$.

The complete eigenstructure of a polynomial matrix satisfies the well-known index sum theorem {\blue (see \cite[Theorem 3]{VergVDK} for the original version for rational matrices, and \cite{DDM} for the more specific polynomial matrix case).}

\begin{theo} \label{thm.indexsum} {\rm (Index Sum Theorem)}
  Let $P(\la)\in \FF[\lambda]^{m\times n}_d$ be a polynomial matrix of normal rank $r$, with invariant polynomials of degrees $\delta_1,\hdots,\delta_r$, with  partial multiplicity sequence at $\infty$ equal to $\gamma_1,\hdots,\gamma_r$, with minimal indices of $\mathcal{N}_\ell (P)$ equal to $\eta_1,\hdots,\eta_{m-r}$ and with minimal indices of $\mathcal{N}_r (P)$ equal to $\varepsilon_1,\hdots,\varepsilon_{n-r}$. Then,
  $$ rd=  \sum_{i=1}^{m-r} \eta_i + \sum_{j=1}^{n-r}\varepsilon_j + \sum_{k=1}^{r} \gamma_k + \sum_{\ell=1}^{r} \delta_\ell .$$
\end{theo}

\subsection{Dual minimal bases and related properties} \label{subsec.dualminbas}
{\blue We now recall the concept of dual minimal bases as defined in \cite{DDMV-zig}, which are closely linked to the classical dual rational subspaces introduced in \cite[Section 6]{For75}.} For brevity, we often say in this paper that a polynomial matrix $M(\la) \in \FF[\la]^{m\times n}$ is a minimal basis if its rows form a minimal basis of the rational subspace they span when $n \geq m$ or if its columns form a minimal basis of the rational subspace they span when $m \geq n$.
\begin{deff}\label{def.dual-m-b}
Two polynomial matrices $M(\lambda)\in \FF[\lambda]^{m\times k}$
and $N(\lambda)\in \FF[\lambda]^{n\times k}$ are dual minimal bases if they are minimal bases satisfying $m+n=k$ and $M(\lambda) \, N(\lambda)^T=0$.
\end{deff}
Observe that the dual minimal bases in Definition \ref{def.dual-m-b} satisfy that the rows of  $M(\lambda)$ form a minimal basis of ${\cal N}_\ell (N(\la)^T)$
and that the columns of $N(\lambda)^T$ form a minimal basis of ${\cal N}_r (M(\la))$. As a consequence, the minimal indices of ${\cal N}_r (M(\la))$ are the degrees of the rows of $N(\la)$ and the minimal indices of ${\cal N}_\ell (N(\la)^T)$ are the degrees of the rows of $M(\la)$.

Dual minimal bases satisfy Theorem \ref{thm.dualbasissum}, whose ``direct part'' was proved in \cite[p. 503]{For75} (see other proofs in \cite[Remark 2.14]{DDMV-zig} and in \cite[Lemma 3.6]{DDVD}) and whose ``converse part'' was proved in \cite[Theorem 6.1]{DDMV-zig}.
\begin{theo} \label{thm.dualbasissum}
  Let $M(\la)\in \FF[\lambda]^{m\times (m+n)}$ and $N(\lambda) \in \FF[\lambda]^{n\times (m+n)}$
  be dual minimal bases with the degrees of their rows equal to $(d_1,\hdots,d_m)$
  and to $(d'_1,\hdots,d'_{n} )$, respectively. Then
\begin{equation} \label{eqn.keyequality}
  \sum_{i=1}^m d_i \,=\, \sum_{j=1}^{n} d'_j \,.
\end{equation}
Conversely, given any two lists of nonnegative integers
$(d_1,\hdots,d_m)$ and $(d'_1,\hdots,d'_{n})$
satisfying \eqref{eqn.keyequality},
there exists a pair of dual minimal bases $M(\la)\in \FF[\lambda]^{m\times (m+n)}$
and $N(\lambda) \in \FF[\lambda]^{n\times (m+n)}$ such that the degrees of the rows of $M(\la)$ and $N(\la)$ are $(d_1,\hdots,d_m)$ and $(d'_1,\hdots,d'_{n})$, respectively.
\end{theo}

A corollary of Theorem \ref{thm.dualbasissum} is the following result.

\begin{cor} \label{cor.leftcolumnmin} Let $P(\la)\in \FF[\lambda]^{m\times n}$ be a polynomial matrix of normal rank $r$, with minimal indices of $\mathcal{N}_\ell (P)$ equal to $\eta_1,\hdots,\eta_{m-r}$, with minimal indices of $\mathcal{N}_r (P)$ equal to $\varepsilon_1,\hdots,\varepsilon_{n-r}$,
with minimal indices of ${\mathcal Row} (P)$ equal to $\rho_1,\ldots, \rho_r$, and with minimal indices of ${\mathcal Col} (P)$ equal to $c_1,\ldots, c_r$. Then
$$
\sum_{i=1}^{m-r} \eta_i = \sum_{i=1}^{r} c_i \qquad \mbox{and} \qquad
\sum_{i=1}^{n-r} \varepsilon_i = \sum_{i=1}^{r} \rho_i \, .
$$
\end{cor}

\begin{proof}
We only prove the first equality, since the second one follows from applying the first to $P(\la)^T$. Let us arrange a minimal basis of $\mathcal{N}_\ell (P)$ as the rows of a matrix $M(\la) \in \FF[\la]^{(m-r) \times m}$ and a minimal basis of ${\mathcal Col} (P)$ as the columns of a matrix $N(\la)^T \in \FF[\la]^{m \times r}$. Then $M(\la) N(\la)^T = 0$, which implies that $M(\la)$ and $N(\la)$ are dual minimal bases and the first equality follows from Theorem \ref{thm.dualbasissum}.
\end{proof}

Combining Theorem \ref{thm.indexsum} and Corollary \ref{cor.leftcolumnmin}, we obtain the following dual version of the Index Sum Theorem.
\begin{cor}
\label{thm.indexsumdual} {\rm (Dual version of the Index Sum Theorem)}
  Let $P(\la)\in \FF[\lambda]^{m\times n}_d$ be a polynomial matrix of normal rank $r$, with invariant polynomials of degrees $\delta_1,\hdots,\delta_r$, with  partial multiplicity sequence at $\infty$ equal to $\gamma_1,\hdots,\gamma_r$, with minimal indices of ${\mathcal Row} (P)$ equal to $\rho_1,\hdots,\rho_r$ and with minimal indices of ${\mathcal Col} (P)$ equal to $c_1,\hdots,c_r$. Then,
  $$ rd=  \sum_{i=1}^{r} c_i + \sum_{j=1}^{r}\rho_j + \sum_{k=1}^{r} \gamma_k + \sum_{\ell=1}^{r} \delta_\ell .$$
\end{cor}

\subsection{Generic complete eigenstructures in $\CC[\la]^{m\times n}_{d,r}$} \label{subsec.genericandrii} We recall in this subsection the main results of \cite{dmy-dop-2017}. For that, we need to introduce some concepts. First, we introduce a distance in the vector space (over the field $\CC$) $\CC[\la]^{m\times n}_{d}$ in terms of the Frobenius matrix norm of complex matrices as follows: Given $P(\lambda) = \lambda^{d}P_{d} + \dots +  \lambda P_1 + P_0 \in \CC[\la]^{m\times n}_{d}$ and $Q(\lambda) = \lambda^{d}Q_{d} + \dots +  \lambda Q_1 + Q_0 \in \CC[\la]^{m\times n}_{d}$, where $P_i, Q_i \in \CC^{m \times n}$, for $i=0, \dots, d$, the distance between $P(\la)$ and $Q(\la)$ is
\begin{equation}\label{ed.defdistance} {\blue
\mathrm{dist} (P,Q)} := \left( \sum_{i=0}^d || P_i - Q_i ||_F^2 \right)^{\frac{1}{2}}.
\end{equation}
This makes $\CC[\la]^{m\times n}_{d}$ a metric space and allows us to define in it limits, open and closed sets, closures of sets and any other topological concept. The closure of any subset $\mathcal{A}$ of $\CC[\la]^{m\times n}_{d}$ will be denoted by $\overline{\mathcal{A}}$.

Given $P(\la) \in \CC[\la]^{m\times n}_{d}$, we define the orbit of $P(\la)$, denoted by $\orb(P)$, as the set of polynomial matrices in $\CC[\la]^{m\times n}_{d}$ with the same complete eigenstructure as $P(\la)$. The closure of $\orb(P)$ is denoted by $\overline{\orb} (P)$. Observe that all the polynomial matrices in $\orb(P)$ have the same rank, since the complete eigenstructure determines the rank, and the same degree, since the first term in the partial multiplicity sequence at $\infty$ determines the degree according to Lemma \ref{lemm.degandinfty}.

The main result in \cite{dmy-dop-2017} describes  $\CC[\la]^{m\times n}_{d,r}$ in terms of closures of orbits of certain polynomial matrices with very particular complete eigenstructures. It is stated in the next theorem {\blue when $r < \min \{m , n \}$}.

\begin{theo} \label{thm.mainthandrii2017} {\rm \cite[Theorem 3.2]{dmy-dop-2017}}
Let $m,n,r$ and $d$ be integers such that $m,n \geq 2$, $d \geq 1$ and $1 \leq r < \min\{m,n\}$. Define $rd + 1$ complete eigenstructures ${\mathbf K}_a$ of  polynomial matrices in $\CC[\la]^{m\times n}_{d,r}$ with $r$ invariant polynomials all equal to one, with all the terms of the partial multiplicity sequence at $\infty$ equal to zero (equivalently, without finite or infinite eigenvalues), with $m-r$ minimal indices of the left null space equal to $\beta$ and $\beta +1$, and with $n-r$ minimal indices of the right null space equal to $\alpha$ and $\alpha +1$, as follows:
\begin{equation}
\label{eq.kcilist}
{\mathbf K}_a : \{\underbrace{\alpha+1, \dots , \alpha+1}_{s},\underbrace{\alpha, \dots , \alpha}_{n-r-s}, \underbrace{\beta+1, \dots , \beta+1}_{t}, \underbrace{\beta, \dots , \beta}_{m-r-t}\}
\end{equation}
for $a = 0,1,\dots,rd$, where $\alpha = \lfloor a / (n-r) \rfloor$, $s = a \mod (n-r)$, $\beta = \lfloor (rd-a)/(m-r) \rfloor$,
and $t = (rd-a) \mod (m-r)$. Then,
\begin{itemize}
\item[\rm (i)] There exists a polynomial matrix $K_a (\la) \in \CC[\la]^{m\times n}_{d,r}$ of degree exactly $d$ and normal rank exactly $r$ with the complete eigenstructure ${\mathbf K}_a$ for $a = 0,1,\dots,rd$;

\item[\rm (ii)]  For every polynomial matrix $M(\la) \in \CC[\la]^{m\times n}_{d,r}$, there exists an integer a such that $\overline{\orb}(K_a)\supseteq \overline{\orb}(M)$;

\item[\rm (iii)] $\overline{\orb}(K_{a}) \bigcap \orb(K_{a'})= \emptyset$ whenever $a\neq a'$;

\item[\rm (iv)]  $\displaystyle \CC[\la]^{m\times n}_{d,r} = \bigcup_{0 \leq a \leq rd} \overline{\orb}(K_a)$ and $\CC[\la]^{m\times n}_{d,r}$ is a closed subset of $\CC[\la]^{m\times n}_{d}$.
\end{itemize}
\end{theo}

Moreover, it was proved in \cite[Corollary 3.3]{dmy-dop-2017} that for each $a =0,1,\ldots , rd$, the orbit $\orb(K_a)$ is an open subset of $\CC[\la]^{m\times n}_{d,r}$ (in the subspace topology of $\CC[\la]^{m\times n}_{d,r}$ corresponding to the distance \eqref{ed.defdistance}). This means that $\bigcup_{0 \leq a \leq rd} \orb(K_a)$ is an open and dense subset of $\CC[\la]^{m\times n}_{d,r}$, which justifies to term the complete eigenstructures in \eqref{eq.kcilist} as the generic eigenstructures of the polynomial matrices in  $\CC[\la]^{m\times n}_{d,r}$. As we have explained in Section \ref{sec.intro}, one of the main objectives of this paper is to provide an alternative description of $\CC[\la]^{m\times n}_{d,r}$ {\blue when $r < \min \{m,n\}$} in terms of the union of the closures of some sets of polynomial matrices that can be factorized in certain specific ways and to relate this description with that in Theorem \ref{thm.mainthandrii2017}-(iv). This is done in Section \ref{sect.compact}.

{\blue Next we consider the generic eigenstructures in the limiting full rank case $r=\min\{m,n\}$, which is not covered by Theorem \ref{thm.mainthandrii2017}. Note that in this case $\CC[\la]^{m\times n}_{d,r} = \CC[\la]^{m\times n}_{d}$ is just the whole set of $m\times n$ polynomial matrices of degree at most $d$. If $m=n$, there is only one generic eigenstructure in $\CC[\la]^{n\times n}_{d}$, which obviously corresponds to regular matrix polynomials, i.e., they do not have minimal indices at all, of degree exactly $d$ and with all their $nd$ eigenvalues distinct. If $m< n$ (resp., $m > n$), there is only one generic eigenstructure which has been described in \cite[Theorem 3.7]{dmy-dop-2017} (resp., \cite[Theorem 3.8]{dmy-dop-2017}). For brevity, we only state the corresponding result for $m<n$, since the result for $m>n$ is analogous and can be obtained by transposition.
\begin{theo} {\rm \cite[Theorem 3.7]{dmy-dop-2017}}
Let $m,n$, $m<n$, and $d$ be positive integers and define the complete eigenstructure ${\mathbf K}_{rp}$ of polynomial matrices in $\CC[\la]^{m\times n}_{d}$ without finite or infinite eigenvalues, without  left minimal indices, and with $n-m$ minimal indices of the right null space equal to $\alpha$ and $\alpha +1$ as follows:
\[
{\mathbf K}_{rp}: \{ \underbrace{\alpha+1, \dots , \alpha+1}_{s},\underbrace{\alpha, \dots , \alpha}_{n-m-s} \} \, ,
\]
where $\alpha = \lfloor md / (n-m) \rfloor$, $s = md \mod (n-m)$. Then,
\begin{itemize}
\item[(i)] There exists a polynomial matrix $K_{rp} \in \CC[\la]^{m\times n}_{d}$ of degree exactly $d$ and normal rank exactly $m$ with the complete eigenstructure ${\mathbf K}_{rp}$;
\item[(ii)] $\CC[\la]^{m\times n}_{d} = \overline{\orb}(K_{rp})$.
\end{itemize}
\end{theo}
Observe that every polynomial in $\orb(K_{rp})$ is, according to  Theorem \ref{thm.minbasischar}, an $m\times n$ minimal basis with the degrees of all its rows equal to $d$, and with the minimal indices of its right null space
differing at most by one. Thus, when $m\ne n$, the polynomial matrices in $\CC[\la]^{m\times n}_{d}$ are generically minimal bases with the degrees of their rows all equal to $d$, if $m<n$, or with the degrees of their columns all equal to $d$, if $m>n$.
}

\subsection{Generic polynomial matrices in $\CC[\la]^{r\times (r+s)}_{\underline{\mathbf{d}}}$ } \label{subsec.robust2summ}
The last subsection in these preliminaries presents a result from \cite{dop-vd-minbasII} that describes the generic polynomial matrices in the vector space (over the field $\CC$) $\CC[\la]^{r\times (r+s)}_{\underline{\mathbf{d}}}$, where $r,s >0$. Observe that if $\underline{\mathbf{d}} = (d_1, d_2, \ldots ,d_r)$ and $d = \max_{1\leq i \leq r} d_i$, then $\CC[\la]^{r\times (r+s)}_{\underline{\mathbf{d}}}$ is a subspace of $\CC[\la]^{r\times (r+s)}_d$ and we can use the distance \eqref{ed.defdistance} in $\CC[\la]^{r\times (r+s)}_{\underline{\mathbf{d}}}$. Moreover, this allows us to define naturally the partial multiplicity sequence at $\infty$ of any $M(\la) \in \CC[\la]^{r\times (r+s)}_{\underline{\mathbf{d}}}$ as the partial multiplicity at $0$ of $\la^d P(1/\la)$.

Next, we define an important subset of $\CC[\la]^{r\times (r+s)}_{\underline{\mathbf{d}}}$, which is proved to be generic in Theorem \ref{thm.genericrobust}.

\begin{deff} \label{def.genericrobust} Let $r,s > 0$ be two positive integers, consider the set $\CC[\la]^{r\times (r+s)}_{\underline{\mathbf{d}}}$, where $\underline{\mathbf{d}} = (d_1, d_2, \ldots ,d_r)$ is a list of nonnegative integers, and define
$$
k' = \left\lceil \frac{\sum_{i=1}^{r} d_i}{s} \right\rceil \quad \mbox{and} \quad
s\, k' = \sum_{i=1}^{r} d_i + t, \quad \mbox{where $0\leq t < s$}.
$$
Then $\mathcal{G}[\la]^{r\times (r+s)}_{\underline{\mathbf{d}}} \subset \CC[\la]^{r\times (r+s)}_{\underline{\mathbf{d}}}$ is the set of polynomial matrices whose $i$th row has degree exactly $d_i$, for $i= 1, \ldots , r$, whose rows form a minimal basis, and such that their right nullspaces have $s$ minimal indices, $t$ of them equal to $k'-1$ and $s-t$ equal to $k'$.
\end{deff}

Theorem 5.3 in \cite{dop-vd-minbasII} implies that $\mathcal{G}[\la]^{r\times (r+s)}_{\underline{\mathbf{d}}}$ is equal to the set of the polynomial matrices that have full-trimmed-Sylvester rank {\blue (the reader can see in \cite[Definition 5.1]{dop-vd-minbasII}, the definition of this concept, though it is not used in this paper)}. Combining this fact with \cite[Theorem 6.2]{dop-vd-minbasII}, we obtain the following result.

\begin{theo} \label{thm.genericrobust} $\mathcal{G}[\la]^{r\times (r+s)}_{\underline{\mathbf{d}}}$ is an open and dense subset of $\CC[\la]^{r\times (r+s)}_{\underline{\mathbf{d}}}$.
\end{theo}

\section{Minimal rank factorizations of polynomial matrices} \label{sec.minfactorizations}
We consider in this section factorizations of a polynomial matrix $P(\la)$ {\it into products of other polynomial matrices} that reveal the normal rank and the degree of $P(\la)$.

\begin{deff} \label{def.rankfactorization} Let $P(\la) \in \FF[\la]^{m\times n}$ with $\rank (P) = r>0$. A factorization of $P(\la)$ as $P(\la) = L(\la) E(\la) R(\la)$ with $L(\la) \in \FF[\la]^{m \times r}, E(\la) \in \FF[\la]^{r \times r}$ and $R(\la) \in \FF[\la]^{r \times n}$ is called a rank factorization of $P(\la)$.
\end{deff}

The name ``rank factorization'' in the definition above reminds us that the sizes of the factors $L(\la), E(\la),$ and $R(\la)$ reveal the rank of $P(\la)$. Standard linear algebra properties of matrices over the field $\FF (\la)$, in particular, the inequality $\rank (L(\la) E(\la) R(\la)) \leq \min \{\rank(L), \rank (E), \rank (R) \}$, immediately imply the following simple well-known results.

\begin{lem} \label{lemm:fact0} Let $P(\la) \in \FF[\la]^{m\times n}$ with $\rank(P) = r>0$ and $P(\la) = L(\la) E(\la) R(\la)$ be a rank factorization of $P(\la)$. Then,
\begin{itemize}
\item[\rm (i)] $\rank (L) = \rank (E) = \rank (R) = r$ and, in particular, $E(\la)$ is nonsingular,

\item[\rm (ii)] ${\mathcal N}_\ell (P) = {\mathcal N}_\ell (L)$,

\item[\rm (iii)] ${\mathcal N}_r (P) = {\mathcal N}_r (R)$,

\item[\rm (iv)] ${\mathcal Row} (P) = {\mathcal Row} (R)$, and the rows of $R(\la)$ are a polynomial basis of ${\mathcal Row} (P)$,

\item[\rm (v)] ${\mathcal Col} (P) = {\mathcal Col} (L)$, and the columns of $L(\la)$ are a polynomial basis of ${\mathcal Col} (P)$.
\end{itemize}
\end{lem}

The following simple lemma is valid for rational matrices {\blue (and, so, for constant and polynomial matrices).}

\begin{lem} \label{lemm:fact0bis} Let $L(\la) \in \FF(\la)^{m \times r}, E(\la) \in \FF(\la)^{r \times r}$ and $R(\la) \in \FF(\la)^{r \times n}$ with $\rank (L) = \rank (R) = r>0$. Then $\rank (L(\la) E(\la) R(\la)) = \rank (E(\la))$.
\end{lem}

{\blue
\begin{proof}
It follows from combining the equalities ${\mathcal N}_r (L(\la) E(\la) R(\la)) = {\mathcal N}_r (E(\la) R(\la))$ and ${\mathcal N}_\ell (E(\la) R(\la)) = {\mathcal N}_\ell (E(\la))$ with the rank-nullity theorem.
\end{proof}
}

{\blue We will often consider rank factorizations with $E (\la) = I_r$. In this case a rank factorization is expressed as the product of just two factors as $P(\la) = L(\la) R(\la)$.}

{\blue Lemma \ref{lemm:fact0}-(iv)-(v) combined with Theorem \ref{thm.mackeystrong} immediately imply the following lower bound for the sum of the degrees of the factors of any rank factorization of $P(\la)$.
\begin{lem} \label{lemm:lowbounddeg}
Let $P(\la) \in \FF[\la]^{m\times n}$ with $\rank (P) = r >0$ and $P(\la) = L(\la) E (\la) R(\la)$ be a rank factorization of $P(\la)$. Let $\rho_{max}$ be the largest minimal index of ${\mathcal Row} (P)$ and $c_{max}$ be the largest minimal index of ${\mathcal Col} (P)$. Then
$$
\deg (L) + \deg (E) + \deg (R) \geq \deg (L) + \deg (R) \geq \rho_{max} + c_{max} .
$$
Thus, if $\rho_{max} + c_{max} > \deg (P)$, then no rank factorization of $P(\la)$, with three or two factors, satisfies that the sum of the degrees of the factors is equal to the degree of $P(\la)$.
\end{lem}
}

An example of a rank factorization of a polynomial matrix can be obtained from truncating the Smith factorization and from elementary properties of matrix multiplication. This is stated in the next lemma.

\begin{lem} \label{lemm:truncsmith} Let $P(\la) \in \FF[\la]^{m\times n}$ with {\blue $\rank(P) = r>0$} and Smith factorization $P(\la) = U(\la) S(\la) V(\la)$ as in \eqref{eq.snf}. Let $L (\la) \in \FF[\la]^{m \times r}$ be the polynomial matrix whose columns are the first $r$ columns of $U(\la)$, $E (\la) \in \FF[\la]^{r \times r}$ be the diagonal polynomial matrix whose diagonal entries are the first $r$ diagonal entries of $S(\la)$, and $R (\la) \in \FF[\la]^{r \times n}$ be the polynomial matrix whose rows are the first $r$ rows of $V(\la)$. Then, $P(\la) = L(\la) E(\la) R(\la)$ is a rank factorization of $P(\la)$.
\end{lem}

\begin{deff} \label{def.Smithrankfactorization} {\blue
A factorization $P(\la) = L(\la) E(\la) R(\la)$ as in Lemma \ref{lemm:truncsmith}} is called a Smith rank factorization of $P(\la)$.
\end{deff}

Smith rank factorizations $P(\la) = L(\la) E(\la) R(\la)$ reveal the invariant polynomials of $P(\la)$ in $E(\la)$, which is a very important information, in addition to the rank of $P(\la)$. However, in general, the columns of $L(\la)$ are not a minimal basis of ${\mathcal Col} (P)$ and the rows of $R(\la)$ are not a minimal basis of ${\mathcal Row} (P)$. Moreover, the degree properties of a Smith rank factorization of $P(\la)$ are not optimal in general. The next example illustrates these facts.

\begin{example} \label{ex.smithrankfact}
{\rm Let
\begin{equation} \label{eq.PexampleSm}
P(\la) = \left[
\begin{array}{ccc}
  \la^8 & 0 & 0 \\
  \la^6 + 1 & -\la^7 & -\la^5 \\
  1 & -\la^7 & -\la^5
\end{array}
\right] \, .
\end{equation}
Then,
\begin{align}
P(\la) & \nonumber =
\left[
\begin{array}{ccc}
  \la^8 & \la^2 & 1 \\
  \la^6 + 1 & 1 & 0 \\
  1 & 0 & 0
\end{array}
\right]
\left[
\begin{array}{ccc}
  1 & 0 & 0 \\
  0 & \la^{11} & 0\\
  0 & 0 & 0
\end{array}
\right]
\left[
\begin{array}{ccc}
  1 & -\la^7 & -\la^5 \\
  0 & \la^2 & 1 \\
  0 & 1 & 0
\end{array}
\right], \\ \label{eq.PexampleFactor}
P(\la) & =
\left[
\begin{array}{cc}
  \la^8 & \la^2  \\
  \la^6 + 1 & 1  \\
  1 & 0
\end{array}
\right]
\left[
\begin{array}{cc}
  1 & 0 \\
  0 & \la^{11}
\end{array}
\right]
\left[
\begin{array}{ccc}
  1 & -\la^7 & -\la^5 \\
  0 & \la^2 & 1
\end{array}
\right] =: L(\la) E(\la) R(\la)
\end{align}
are, respectively, a Smith factorization and a Smith rank factorization of $P(\la)$. Note that according to Theorem \ref{thm.minbasischar} neither $L(\la)$ nor $R(\la)$ are minimal bases because their highest-column-degree and highest-row-degree coefficients are, respectively,
$$
L_{hc} = \begin{bmatrix}
           1 & 1 \\
           0 & 0 \\
           0 & 0
         \end{bmatrix}
\quad \mbox{and} \quad
R_{hr} = \begin{bmatrix}
           0 & -1 & 0 \\
           0 & 1 & 0
         \end{bmatrix},
$$
which do not have full column and full row rank, respectively.
Observe that the degree of $P(\la)$ is $8$ and that is not equal to the sum of the degrees of the three factors in \eqref{eq.PexampleFactor}, which is $26$. This inequality is expected because the entries with highest degrees in each factor do not interact when the product $L(\la) E(\la) R(\la)$ is computed. But note also that if  \eqref{eq.PexampleFactor} is expanded into a sum of rank one matrices as follows
\begin{equation} \label{eq.PexampleFactor-EXP}
P(\la)  =
\left[
\begin{array}{c}
  \la^8  \\
  \la^6 + 1   \\
  1
\end{array}
\right]
\left[
\begin{array}{c}
  1
\end{array}
\right]
\left[
\begin{array}{ccc}
  1 & -\la^7 & -\la^5
\end{array}
\right] +
\left[
\begin{array}{c}
\la^2  \\
1  \\
0
\end{array}
\right]
\left[
\begin{array}{c}
  \la^{11}
\end{array}
\right]
\left[
\begin{array}{ccc}
  0 & \la^2 & 1
\end{array}
\right],
\end{equation}
then the degrees of both terms are equal to $15$, again much larger than $\deg(P) = 8$.
{\blue \qed}}
\end{example}

In the rest of this section, we explore other rank factorizations, different from Smith rank factorizations, of a polynomial matrix $P(\la)$ whose factors provide minimal bases of ${\mathcal Col} (P)$ and/or ${\mathcal Row} (P)$ and reveal the degree of $P(\la)$. We emphasize that, in general, such factorizations do not reveal explicitly the invariant polynomials of $P(\la)$.

We will need in the sequel the two auxiliary Lemmas \ref{lemm.degreerank1} and \ref{lemm:fact1}. Lemma \ref{lemm.degreerank1} implies, in particular, that any rank factorization of a polynomial matrix $P(\la)$ with normal rank equal to one reveals the degree of $P(\la)$ via the sum of the degrees of the three factors. The simple proof of this lemma is omitted.

\begin{lem} \label{lemm.degreerank1}
 Let $L(\la)\in \FF[\la]^{m\times 1}, E(\la)\in \FF[\la]^{1\times 1}, R(\la)\in \FF[\la]^{1\times n}$ and $P(\la)=L(\la)E(\la)R(\la)$. Then
$\deg (P) = \deg(L) + \deg (E) + \deg(R).$
\end{lem}

Lemma \ref{lemm:fact1} is a consequence of \cite[Theorem 2.5.7]{wolovich}, which introduces an algorithm for transforming any polynomial matrix with full column rank into a column reduced polynomial matrix via multiplication on the right by unimodular matrices. In order to be self-contained, we include a short proof of this lemma.

\begin{lem} \label{lemm:fact1} \phantom{llpp}
\begin{itemize}
\item[\rm (i)] Let $L(\la)\in \FF[\la]^{m\times r}$ be a polynomial matrix such that the constant matrix $L(\la_0)$ has full column rank $r$ for all $\la_0\in \overline{\FF}$. Then, $L(\la)$ can be factorized as $L(\la)=L_c(\la) V(\la)$, where the columns of $L_c(\la)\in \FF[\la]^{m\times r}$ form a minimal basis of ${\mathcal Col} (L)$ and  $V(\la)\in \FF[\la]^{r\times r}$ is unimodular. Hence, the degrees of the columns of $L_c (\la)$ are the minimal indices of ${\mathcal Col} (L)$.

\item[\rm (ii)] Let $R(\la)\in \FF[\la]^{r\times n}$ be a polynomial matrix such that the constant matrix $R(\la_0)$ has full row rank $r$ for all $\la_0\in \overline{\FF}$. Then, $R(\la)$ can be factorized as $R(\la)=U(\la)R_r(\la)$, where the rows of $R_r(\la) \in \FF[\la]^{r\times n}$ form a minimal basis of  ${\mathcal Row} (R)$ and  $U(\la)\in \FF[\la]^{r\times r}$ is unimodular. Hence, the degrees of the rows of $R_r(\la)$ are the minimal indices of ${\mathcal Row} (R)$.
\end{itemize}
\end{lem}

\begin{proof} We only prove item (i), since item (ii) is obtained from item (i) by transposition. The columns of $L(\la)$ are a basis of ${\mathcal Col} (L)$. If the columns of $L_c (\la)$ are any minimal basis of ${\mathcal Col} (L)$, then $L(\la)=L_c(\la) V(\la)$, with $V(\la)$ an $r \times r$ polynomial matrix according to \cite[p. 495]{For75}. In addition, $V(\la)$ must be unimodular since, otherwise,  $L(\la_0)=L_c(\la_0) V(\la_0)$ would have rank strictly smaller than $r$ for any root $\la_0$ of $\det V(\la_0)$.
\end{proof}

The next example illustrates Lemma \ref{lemm:fact1}.
\begin{example} \label{ex.minimaluni}
{\rm The matrices $L(\lambda)$ and $R(\lambda)$ in \eqref{eq.PexampleFactor} can be factorized as follows:
\begin{align} \label{ex.factorL}
\left[
\begin{array}{cc}
  \la^8 & \la^2  \\
  \la^6 + 1 & 1  \\
  1 & 0
\end{array}
\right] & =
\left[
\begin{array}{cc}
  0 & \la^2  \\
  1 & 1  \\
  1 & 0
\end{array}
\right]
\left[
\begin{array}{cc}
  1 & 0  \\
  \la^6 & 1
\end{array}
\right] =: L_c (\la) V(\la),
\\ \label{ex.factorU}
\left[
\begin{array}{ccc}
  1 & -\la^7 & -\la^5 \\
  0 & \la^2 & 1
\end{array}
\right] & =
\left[
\begin{array}{ccc}
  1 & -\la^5 \\
  0 &  1
\end{array}
\right]
\left[
\begin{array}{ccc}
  1 & 0 & 0 \\
  0 & \la^2 & 1
\end{array}
\right] =: U(\la) R_r(\la).
\end{align}
Theorem \ref{thm.minbasischar} implies that the columns of $L_c(\la)$ are a minimal basis, as well as the rows of $R_r(\la)$. Obviously $V(\la)$ and $U(\la)$ are unimodular matrices.
{\blue \qed}}
\end{example}

Theorem \ref{thm.minfact0} presents for each polynomial matrix three different types of rank factorizations, with $E(\la)= I_r$ in the case of items (ii) and (iii).

\begin{theo} \label{thm.minfact0}
Let $P(\la)\in \FF[\la]^{m\times n}$ with {\blue $\rank (P) = r>0$.} Then, $P(\la)$ can be factorized as follows:
\begin{itemize}
\item[\rm (i)] $P(\la)=L_c (\la) E(\la) R_r(\la)$, where $L_c(\la)\in \FF[\la]^{m\times r}$, $E(\la)\in \FF[\la]^{r\times r}$, $R_r(\la)\in \FF[\la]^{r\times n}$, the columns of $L_c (\la)$ form a minimal basis of ${\mathcal Col} (P)$, the rows of $R_r (\la)$ form a minimal basis of ${\mathcal Row} (P)$, and the invariant polynomials of $E(\la)$ are the invariant polynomials of $P(\la)$.

\item[\rm (ii)] $P(\la)=L_c (\la) R(\la)$, where $L_c(\la)\in \FF[\la]^{m\times r}$, $R (\la)\in \FF[\la]^{r\times n}$, the columns of $L_c (\la)$ form a minimal basis of ${\mathcal Col} (P)$ and the invariant polynomials of $R(\la)$ are the invariant polynomials of $P(\la)$.

\item[\rm (iii)] $P(\la)=L(\la) R_r(\la)$, where $L(\la)\in \FF[\la]^{m\times r}$, $R_r(\la)\in \FF[\la]^{r\times n}$, the rows of $R_r (\la)$ form a minimal basis of ${\mathcal Row} (P)$, and the invariant polynomials of $L(\la)$ are the invariant polynomials of $P(\la)$.
\end{itemize}
\end{theo}
\begin{proof} Let $P(\la)=\widetilde{L} (\la) \widetilde{E}(\la) \widetilde{R}(\la)$ with $\widetilde{L}(\la)\in \FF[\la]^{m\times r}$, $\widetilde{E}(\la)\in \FF[\la]^{r\times r}$, and $\widetilde{R}(\la)\in \FF[\la]^{r\times n}$, be a Smith rank factorization as in Lemma \ref{lemm:truncsmith}. Therefore, $\widetilde{L} (\la_0)$ and $\widetilde{R}(\la_0)$ have, respectively, full column rank and full row rank for all $\la_0 \in \overline{\FF}$, because they are formed by columns and rows, respectively, of unimodular matrices. Then using the factorizations in Lemma \ref{lemm:fact1} applied to $\widetilde{L} (\la)$ and $\widetilde{R}(\la)$, we get the following three expressions,
\begin{align}
  P(\la) & = L_c (\la) \, (V(\la) \widetilde{E} (\la) U(\la)) \,  R_r (\la), \label{eq.auxxmin1}\\
  P(\la) & = L_c (\la) \, (V(\la) \widetilde{E}(\la) \widetilde{R}(\la)) ,  \label{eq.auxxmin2} \\
  P(\la) & = (\widetilde{L}(\la) \widetilde{E}(\la) U(\la) ) \, R_r( \la) .  \label{eq.auxxmin3}
\end{align}
The factorization in \eqref{eq.auxxmin1} proves item (i) with $E(\la) = V(\la) \widetilde{E} (\la) U(\la)$, because the $r$ diagonal entries of $\widetilde{E} (\la)$ are the invariant polynomials of $P(\la)$ and they do not change under multiplications by unimodular matrices. The statements about ${\mathcal Col} (P)$ and ${\mathcal Row} (P)$ follow from Lemma \ref{lemm:fact0}.

The factorization in \eqref{eq.auxxmin2} proves item (ii) with $R (\la) = V(\la) \widetilde{E}(\la) \widetilde{R}(\la)$. Note that $\widetilde{R}(\la)$ is formed by the first $r$ rows of a unimodular matrix $\widetilde{V}(\la) \in \FF[\la]^{n\times n}$, according to Lemma \ref{lemm:truncsmith}. Thus,
$$R (\la) = V(\la) \begin{bmatrix}
                                 \widetilde{E}(\la) & 0
                               \end{bmatrix} \widetilde{V}(\la),$$
and indeed the invariant polynomials of $R(\la)$ are the same of those of $\widetilde{E}(\la)$, which in turn are those of $P(\la)$. The statement about ${\mathcal Col} (P)$ follows again from Lemma \ref{lemm:fact0}.

Analogously, the factorization in \eqref{eq.auxxmin3} proves item (iii).
\end{proof}

\begin{deff} \label{def.minrankfactorization} Any of the three factorizations introduced in Theorem \ref{thm.minfact0} is called a minimal rank factorization of $P(\la)$.
\end{deff}

The name ``minimal rank factorization'' in Definition \ref{def.minrankfactorization} reminds us that these factorizations display a minimal basis of ${\mathcal Col} (P)$ and/or a minimal basis of ${\mathcal Row} (P)$, in addition to the rank of $P(\la)$.

\begin{rem} \label{rem.uniqueminfact} {\rm The minimal rank factorizations in Theorem \ref{thm.minfact0} are not unique. In fact  $L_c (\la)$ can be any of {\blue all possible} minimal bases of ${\mathcal Col} (P)$ and $R_r (\la)$ can be any of {\blue all possible} minimal bases of ${\mathcal Row} (P)$. However, note that once $L_c (\la)$ and/or $R_r (\la)$ are chosen, $E(\la)$ in item (i) is uniquely determined by this choice, $R(\la)$ in item (ii) is uniquely determined by this choice, and $L(\la)$ in item (iii) is uniquely determined by this choice. {\blue This follows from the fact that $L_c(\la)$ has a polynomial left inverse and $R_r (\la)$ has a polynomial right inverse.}}
\end{rem}

{\blue
\begin{rem} {\rm (Numerical computation of minimal rank factorizations) If $\FF = \mathbb{R}$ or $\mathbb{C}$, minimal rank factorizations of a polynomial matrix $P(\la)$ can be computed by using any of the existing algorithms that compute efficiently minimal bases of the null spaces of a polynomial matrix by applying unitary transformations on constant matrices. There are different families of such algorithms. For instance, some are based on Sylvester resultant matrices \cite{antoniouLAA2005} and others are based on computing minimal bases of the null spaces of a linearization of $P(\la)$ via the staircase algorithm \cite{beelen1987,nofervdlaa2023,vd1979} and recovering the minimal bases of the null spaces of $P(\la)$ from those of the linearization \cite{fdtfdmackeyELA}. Then, a minimal basis $L_c (\la)$ of ${\mathcal Col} (P)$ can be computed in two steps as follows: (1) Compute a minimal basis $Q(\la)$ of ${\mathcal N}_\ell (P)$; (2) Compute a minimal basis $L_c (\la)$ of ${\mathcal N}_r (Q) = {\mathcal Col} (P)$. Then, $R (\la)$ in Theorem \ref{thm.minfact0}-(ii) can be computed by solving a system of linear equations with unique solution for its coefficients.
Applying this procedure to $R(\la)^T$ yields the factorization in Theorem \ref{thm.minfact0}-(i). The factorization in Theorem \ref{thm.minfact0}-(iii) can be obtained in a similar manner. Finally, we mention that the algorithm summarized in \cite[Theorem 4.3]{nofervdsimax2023} computes directly a factorization that is ``almost'' the same as the one in Theorem \ref{thm.minfact0}-(ii) using a unitary decomposition of a generalized state-space model of $P(\la)$. The only missing property is that the computed $L_c (\la)$ is not guaranteed to be column reduced. It remains as an open problem to adapt the algorithm in \cite{nofervdsimax2023} to ensure the column reduced property.
}
\end{rem}
}

{\blue \begin{rem} \label{rem.fullrank} {\rm (Minimal rank factorizations of full rank polynomial matrices) If $P(\la) \in \FF [\la]^{m\times n}$ has full rank $r = \min \{m,n\}$, the minimal rank factorizations in Theorem \ref{thm.minfact0} are simpler. Note that if $r=m$ (resp., $r=n$), then ${\mathcal Col}(P) = \FF(\la)^{m \times 1}$ (resp., ${\mathcal Row}(P) = \FF(\la)^{1 \times n}$). Therefore, if $r=m$, then the minimal basis $L_c (\la) \in \FF^{m\times m}$ in Theorem \ref{thm.minfact0} must be a constant $m\times m$ invertible matrix, that may be taken equal to any of such matrices. In particular, one can take $L_c (\la) = I_m$. Analogously,  if $r=n$, then the minimal basis $R_r (\la) \in \FF^{n\times n}$ in Theorem \ref{thm.minfact0} must be a constant $n\times n$ invertible matrix, that may be taken equal to any of such matrices. In particular, one can take $R_r (\la) = I_n$.}
\end{rem}}

The next example illustrates Theorem \ref{thm.minfact0}.

\begin{example} \label{ex.minimalrankfact}
{\rm In this example, the Smith rank factorization in \eqref{eq.PexampleFactor} is combined with the factorizations in \eqref{ex.factorL} and \eqref{ex.factorU} to obtain the following minimal rank factorizations of $P(\la)$ in \eqref{eq.PexampleSm}:
\begin{align}\label{eq.exminfact1}
  P(\la) & =
\left[
\begin{array}{cc}
  0 & \la^2  \\
  1 & 1  \\
  1 & 0
\end{array}
\right]
\left[
\begin{array}{cc}
  1 & -\la^5  \\
  \la^6 & 0
\end{array}
\right]
\left[
\begin{array}{ccc}
  1 & 0 & 0 \\
  0 & \la^2 & 1
\end{array}
\right]
=: L_c (\la) F(\la) R_r (\la),\\ \label{eq.exminfact2}
P(\la)  &= \left[
\begin{array}{cc}
  0 & \la^2  \\
  1 & 1  \\
  1 & 0
\end{array}
\right]
\left[
\begin{array}{ccc}
  1 & -\la^7 & -\la^5 \\
  \la^6 & 0 & 0
\end{array}
\right]
=: L_c (\la) R(\la), \\ \label{eq.exminfact3}
P(\la)  & =
\left[
\begin{array}{cc}
  \la^8 & 0 \\
  \la^6 + 1 & -\la^5  \\
  1 & -\la^5
\end{array}
\right]
\left[
\begin{array}{ccc}
  1 & 0 & 0 \\
  0 & \la^2 & 1
\end{array}
\right]
=: L(\la) R_r (\la).
\end{align}
The factorizations in \eqref{eq.exminfact1}, \eqref{eq.exminfact2} and \eqref{eq.exminfact3} illustrate, respectively, items (i), (ii) and (iii) of Theorem \ref{thm.minfact0}. Observe that none of them reveals by inspection the invariant polynomials $1$ and $\la^{11}$ of $P(\la)$ in \eqref{eq.PexampleSm}. However, the degree of $P(\la)$, which is $8$, is revealed as the largest degree of the terms in each of the expansions of $P(\la)$ into a sum of rank one matrices stemming from \eqref{eq.exminfact1}, \eqref{eq.exminfact2} and \eqref{eq.exminfact3}. These expansions are the following ones:
\begin{align*}
  P(\la)  & =  \left[
\begin{array}{c}
  0   \\
  1  \\
  1
\end{array}
\right]
\left[
\begin{array}{c}
  1
\end{array}
\right]
\left[
\begin{array}{ccc}
  1 & 0 & 0
\end{array}
\right] + \left[
\begin{array}{c}
  0  \\
  1  \\
  1
\end{array}
\right]
\left[
\begin{array}{c}
  -\la^5
\end{array}
\right]
\left[
\begin{array}{ccc}
 0 & \la^2 & 1
\end{array}
\right] + \left[
\begin{array}{c}
 \la^2  \\
 1  \\
 0
\end{array}
\right]
\left[
\begin{array}{c}
  \la^6
\end{array}
\right]
\left[
\begin{array}{ccc}
  1 & 0 & 0
\end{array}
\right] \\
  & = \left[
\begin{array}{c}
  0  \\
  1  \\
  1
\end{array}
\right]
\left[
\begin{array}{ccc}
  1 & -\la^7 & -\la^5
\end{array}
\right] +
\left[
\begin{array}{c}
\la^2  \\
 1  \\
 0
\end{array}
\right]
\left[
\begin{array}{ccc}
\la^6 & 0 & 0
\end{array}
\right]\\
  & =
\left[
\begin{array}{c}
  \la^8 \\
  \la^6 + 1   \\
  1
\end{array}
\right]
\left[
\begin{array}{ccc}
  1 & 0 & 0
\end{array}
\right]
+
\left[
\begin{array}{c}
  0 \\ -\la^5  \\ -\la^5
\end{array}
\right]
\left[
\begin{array}{ccc}
  0 & \la^2 & 1
\end{array}
\right].
\end{align*}
The term with highest degree in each of these expansions has degree $8$, which is precisely the degree of the polynomial. This behaviour is in contrast with the degrees of the terms in the expansion \eqref{eq.PexampleFactor-EXP} coming from the Smith rank factorization \eqref{eq.PexampleFactor}. This result about degrees holds for any minimal rank factorization and will be proved in {\blue Corollary \ref{thm.degreepredic}.} {\blue \qed}}
\end{example}


{\blue Example \ref{ex.minimalrankfact} motivates us to state, in Corollary \ref{thm.degreepredic}, some degree properties of products of two and three polynomial matrices which are direct consequences of the classical ``predictable-degree property'' for column reduced and row reduced matrices, and, so, for minimal bases \cite[Condition 4(b), p. 495 and Remark 3, p. 497]{For75} and \cite[Thm. 6.3-13, p. 387]{Kai80}.  Recall the following notation introduced in Section \ref{sec.prelim}: $X_{*i}$ denotes the $i$th column of the matrix $X$ and $Y_{i*}$ denotes the $i$th row of $Y$.}

{\blue
\begin{cor}{\rm (Predictable-degree properties for matrix products).} \label{thm.degreepredic}
\begin{itemize}
  \item[\rm (i)] Let $P(\la)=L(\la)R(\la)$, where $L(\la)\in \FF[\la]^{m\times r}$ and $R(\la)\in \FF[\la]^{r\times n}$. If $L(\la)$ is column reduced or $R(\la)$ is row reduced, then
$$ \deg (P) = \max_{1\leq i \leq r} \{ \deg(L_{*i}) + \deg(R_{i*}) \}.$$
  \item[\rm (ii)] Let $ P(\la)=L(\la)E(\la)R(\la)$, where $L(\la)\in \FF[\la]^{m\times r}$, $E(\la)\in \FF[\la]^{r\times s}$ and $R(\la)\in \FF[\la]^{s\times n}$. If $L(\la)$ is column reduced and $R(\la)$ is row reduced, then
$$\deg (P) = \max_{ {\footnotesize \begin{array}{c}  1\leq i \leq r  \\ 1\leq j \leq s \end{array}} } \{\deg(L_{*i}) + \deg (e_{ij}) + \deg(R_{j*}) \}.$$
\end{itemize}
\end{cor}

\begin{proof}
Proof of item (i). Consider first that $L(\la)$ is column reduced. Then, the classical predictable degree property in \cite[Thm. 6.3-13, p. 387]{Kai80} ensures that $\deg (P_{*j}) = \max_{1\leq i \leq r} \{ \deg(L_{*i}) + \deg(R_{ij}) \}$, from which the result follows by taking the maximum over $1\leq j \leq n$. If $R(\la)$ is row reduced, the sought equality is deduced from applying to $P(\la)^T$ the result just proved when the first factor is column reduced.

\medskip

Proof of item (ii). Define $Q(\la) := L(\la) E(\la)$ and express $P(\la) = Q(\la) R(\la)$. Since $R(\la)$ is row reduced, item (i) implies that $\deg (P) = \max_{1\leq j \leq s} \{ \deg(Q_{*j}) + \deg(R_{j*}) \}$. Moreover, since $L(\la)$ is column reduced, item (i) also implies $\deg (Q_{*j}) = \max_{1\leq i \leq r} \{ \deg(L_{*i}) + \deg(e_{ij}) \}$. The result in item (ii) follows from
combining both equalities.
\end{proof}
}

\begin{rem} {\rm Observe that $P(\la)$ in {\blue item (i) of Corollary \ref{thm.degreepredic}} can be expanded as a sum of rank one polynomial matrices as $P(\la) = \sum_{i=1}^{r} L_{*i} (\la) \, R_{i*} (\la)$, while $P(\la)$ in {\blue item (ii)} can be expanded as a sum of rank one polynomial matrices as $P(\la) = \sum_{i=1}^{r} \sum_{j=1}^{s} L_{*i} (\la) \, e_{ij} (\la) \, R_{j*} (\la)$. Thus, taking into account Lemma \ref{lemm.degreerank1}, {\blue Corollary \ref{thm.degreepredic}} states that the degree of $P(\la)$ is precisely the degree of the term(s) with highest degree in such expansions.
}
\end{rem}

In the last part of this section, we study minimal rank factorizations of polynomial matrices that have no finite or infinite eigenvalues. The motivation for this study comes from Theorem \ref{thm.mainthandrii2017}, which shows that rank deficient polynomial matrices have no finite or infinite eigenvalues, generically, when $\FF = \CC$. {\blue We will see that the minimal rank factorizations have stronger properties in this case.}
\begin{theo} \label{thm.minfact} Let $P(\la) \in \FF[\la]^{m\times n}_d$ and $r$ be an integer such that {\blue $0 < r \leq \min\{m,n\}$}. $P(\la)$ has normal rank $r$, degree exactly $d$, and has no finite or infinite eigenvalues
if and only if $P(\la)$ can be factorized as
\begin{equation} \label{eq.minfact1}
P(\la)=L (\la)R (\la), \quad L(\la)\in \FF[\la]^{m\times r},  \quad R(\la)\in \FF[\la]^{r\times n},
\end{equation}
where the columns of $L (\la)$ are a minimal basis, the rows of $R (\la)$ are a minimal basis, and
\begin{equation} \label{eq.minfact2}
\deg (L_{*i}) + \deg (R_{i*} ) =d , \quad \mbox{for $i=1,\hdots,r$.}
\end{equation}
\end{theo}
\begin{proof}
Sufficiency. If $P(\la)$ satisfies \eqref{eq.minfact1} and \eqref{eq.minfact2} with $L (\la)$ and $R(\la)$ minimal bases, then $\rank(P) =r$ follows from Lemma \ref{lemm:fact0bis} with $E(\la) = I_r$, and $\deg(P) =d$ follows from {\blue Corollary \ref{thm.degreepredic}.} Lemma \ref{lemm:fact0} implies that the degrees of the columns of $L(\la)$ are the minimal indices of ${\mathcal Col} (P)$ and that the degrees of the rows of $R(\la)$ are the minimal indices of ${\mathcal Row} (P)$, and \eqref{eq.minfact2} implies that the sum of all these minimal indices is equal to $r d$. Combining this fact with Corollary \ref{thm.indexsumdual}, we see that all the terms of the partial multiplicity sequence at $\infty$ of $P(\la)$ must be zero, as well as all the degrees of the invariant polynomials of $P(\la)$. This is equivalent to state that $P(\la)$ has no finite or infinite eigenvalues.

Necessity. If $P(\la)$ has normal rank $r$, degree $d$ and has no finite or infinite eigenvalues, then we start from a minimal rank factorization of $P(\la)$ given by Theorem \ref{thm.minfact0}-(ii). That is, $P(\la)=L(\la) \widehat{R}(\la)$, where the columns of $L(\la)\in \FF[\lambda]^{m \times r}$ form a minimal basis of ${\mathcal Col} (P)$ and the degrees of these columns, denoted by $c_1,\ldots ,c_r$, are the minimal indices of  ${\mathcal Col} (P)$. Note also that Lemma \ref{lemm:fact0} guarantees that the rows of $\widehat{R}(\la) \in \FF[\lambda]^{r \times n}$ form a basis of ${\mathcal Row} (P)$. Let $\hat \rho_1, \ldots , \hat \rho_r$ be the degrees of the rows of $\widehat{R}(\la)$, which are not necessarily the minimal indices ${\mathcal Row} (P)$. Therefore, their sum is larger than or equal to the sum of the minimal indices $\rho_1, \ldots, \rho_r$ of ${\mathcal Row} (P)$ by the definition of minimal basis. That is
\begin{equation} \label{rhat}  \sum_{i=1}^r \hat \rho_i \ge  \sum_{i=1}^r \rho_i
\end{equation}
and, simultaneously, from {\blue Corollary \ref{thm.degreepredic}-(i)},
\begin{equation} \label{eq.dgeqrhat}
d \geq c_i + \hat \rho_i \qquad \mbox{for $i=1,\ldots, r$}.
\end{equation}
On the other hand Corollary \ref{thm.indexsumdual} implies
\begin{equation} \label{eq.dualistred}
\sum_{i=1}^{r} c_i + \sum_{j=1}^{r}\rho_j  = rd,
\end{equation}
since $P(\la)$ has no finite or infinite eigenvalues, which is equivalent to
$\sum_{k=1}^{r} \gamma_k + \sum_{\ell =1}^{r}\delta_\ell = 0$.
The combination of \eqref{rhat}, \eqref{eq.dgeqrhat} and \eqref{eq.dualistred} leads to
\begin{equation} \label{eq.finalfactnovals}
\sum_{i=1}^{r} c_i + \sum_{j=1}^{r} \hat \rho_j \geq \sum_{i=1}^{r} c_i + \sum_{j=1}^{r}\rho_j  = rd \geq \sum_{i=1}^{r} c_i + \sum_{j=1}^{r} \hat \rho_j \, .
\end{equation}
Therefore, $\sum_{i=1}^r \hat \rho_i =  \sum_{i=1}^r \rho_i$, which implies that $\widehat{R}(\la)$ is a minimal basis of ${\mathcal Row} (P)$. Moreover, \eqref{eq.finalfactnovals} implies $\sum_{i=1}^{r} c_i + \sum_{j=1}^{r} \hat \rho_j = rd$, which combined with \eqref{eq.dgeqrhat} yields
$$ c_i+\hat \rho_i=d, \quad i=1,\hdots, r \, .$$
This completes the proof.
\end{proof}

\begin{rem} {\rm Observe that the hypothesis that {\blue $P(\la) \in \FF [\la]^{m\times n}_d$} has degree exactly $d$ in Theorem \ref{thm.minfact} is redundant, because Lemma \ref{lemm.degandinfty} combined with the hypothesis that $P(\la)$ has not eigenvalues at $\infty$ implies that $\deg(P) =d$. We have included this redundant hypothesis for emphasizing this key property of the polynomial matrices satisfying Theorem \ref{thm.minfact}.}
\end{rem}

We remark that the proof of the {\it necessity} in Theorem \ref{thm.minfact} proves, in fact, that {\it for any} minimal rank factorization as in Theorem \ref{thm.minfact0}-(ii) of a polynomial matrix $P(\la) \in \FF[\la]^{m\times n}_d$ with normal rank $r$, with degree exactly $d$, and without finite or infinite eigenvalues, the factor $R(\la)$ must be a minimal basis and that the degree constraints \eqref{eq.minfact2} must be satisfied. A complementary result can be proved {\it for any} minimal rank factorization as in Theorem \ref{thm.minfact0}-(iii) just by transposing the argument above. These discussions can be formalized into the following theorem.

\begin{theo} \label{thm.corminfact} Let $P(\la) \in \FF[\la]^{m\times n}_d$ be a polynomial matrix with normal rank {\blue $r>0$,} with degree exactly $d$, and without eigenvalues, finite or infinite. Then, the following statements hold:
\begin{itemize}
  \item[\rm (i)] If the minimal rank factorization $P(\la) = L(\la) E(\la) R(\la)$ satisfies the properties in Theorem \ref{thm.minfact0}-(i), then the rows of $\widehat{R} (\la) = E(\la) R(\la)$ form a minimal basis of ${\mathcal Row} (P)$ and the columns of $\widehat{L} (\la) =L(\la)E(\la)$ form a minimal basis of ${\mathcal Col} (P)$. Moreover,
      $$
      \deg (L_{*i}) + \deg (\widehat{R}_{i*}) =  \deg (\widehat{L}_{*i}) + \deg (R_{i*}) = d, \quad \mbox{for $i=1,\ldots r$.}
      $$
  \item[\rm (ii)] If the minimal rank factorization $P(\la) = L(\la) R(\la)$ satisfies the properties in Theorem \ref{thm.minfact0}-(ii), then the rows of $R(\la)$ form a minimal basis of ${\mathcal Row} (P)$. Moreover,
      $$
      \deg (L_{*i}) + \deg (R_{i*}) =  d, \quad \mbox{for $i=1,\ldots r$.}
      $$
  \item[\rm (iii)] If the minimal rank factorization $P(\la) = L(\la) R(\la)$ satisfies the properties in Theorem \ref{thm.minfact0}-(iii), then the columns of $L(\la)$ form a minimal basis of ${\mathcal Col} (P)$. Moreover,
      $$
      \deg (L_{*i}) + \deg (R_{i*}) =  d, \quad \mbox{for $i=1,\ldots r$.}
      $$
\end{itemize}
\end{theo}

The next example illustrates that polynomial matrices without finite nor infinite eigenvalues have rank factorizations that are not minimal, that do not satisfy the degree conditions \eqref{eq.minfact2} and whose factors can have arbitrarily large degrees. Thus, for polynomial matrices without eigenvalues, minimal rank factorizations are clearly preferable.

\begin{example} \label{ex.arbitrarydeg} {\rm Consider the following polynomial matrix $P(\la) \in \CC[\la]^{3 \times 3}_{6}$ with $\rank(P) = 2$ and its following factorizations:
\begin{align}\label{eq.exdeg1}
  P(\la) =
\begin{bmatrix}
  \la^6 & \la^5 & 0 \\
  \la & \la^6 + 1 & \la^2 \\
  0 & \la^4 & 1
\end{bmatrix}
  & =
\left[
\begin{array}{cc}
  \la^5 & 0  \\
  1 & \la^2  \\
  0 & 1
\end{array}
\right]
\left[
\begin{array}{ccc}
  \la & 1 & 0 \\
  0 & \la^4 & 1
\end{array}
\right],\\ \label{eq.exdeg2}
&= \left[
\begin{array}{cc}
  \la^5 & 0  \\
  \la^{p+2} + 1& \la^2  \\
  \la^p & 1
\end{array}
\right]
\left[
\begin{array}{ccc}
  \la & 1 & 0 \\
  -\la^{p+1} & -\la^p + \la^4 & 1
\end{array}
\right],
\end{align}
where $p>3$ is an integer. Both factors in the factorization \eqref{eq.exdeg1} are minimal bases, according to Theorem \ref{thm.minbasischar}, and they satisfy \eqref{eq.minfact2}. This proves that $P(\la)$ has no finite or infinite eigenvalues. In contrast, the left and right factors in \eqref{eq.exdeg2} are not, respectively, column and row reduced polynomial matrices. So, they are not minimal bases. Nevertheless, the factorization in \eqref{eq.exdeg2} is a rank factorization of $P(\la)$. Its factors have arbitrarily high degrees for arbitrarily large values of $p$.
{\blue Lemma \ref{lemm:lowbounddeg} applied to the polynomial matrix $P(\la)$ in this example yields a lower bound $\rho_{max} + c_{max} = 9 > 6 = \deg (P)$, which implies that no rank factorization of $P(\la)$ satisfies that
the sum of the degrees of the factors is equal to the degree of $P(\la)$.}
{\blue \qed}
}
\end{example}

By Theorem \ref{thm.mainthandrii2017}, rank deficient complex polynomial matrices have, generically, no finite or infinite eigenvalues. Combining this with Theorem \ref{thm.minfact}, we obtain that, generically, rank deficient complex polynomial matrices have minimal rank factorizations as simple as those appearing in Theorem \ref{thm.minfact}. However, the degree condition \eqref{eq.minfact2} still allows for a lot of freedom for the possible degrees of the columns of $L(\la)$ and the rows of $R(\la)$, {\blue except when $d$ and $r$ are very small. More precisely, taking into account that the degrees of the columns of $L(\la)$ are determined by the degrees of the rows of $R(\la)$ and that the order of the columns of $L(\la)$ and the rows of $R(\la)$ does not affect the product $L(\la) R(\la)$, the number of different degree distributions in \eqref{eq.minfact1} is equal to the number of $r$-combinations with repetitions from the set of $d+1$ possible values of the degree of each row of $R(\la)$. This amounts to $\binom{d+r}{r}$ different degree distributions. This number is huge except for very small values of $d$ and $r$ and makes it unfeasible to develop a set of parametrizations for the generic set of $m\times n$ polynomial matrices of degree at most $d$ without eigenvalues. However, in the next section, we prove that, generically, the degrees are considerably more constrained and that there are only $rd+1$ different ways to distribute the degrees, all of them characterized by the fact that the degrees of the columns of $L(\la)$ differ at most by one and the degrees of the rows of $R(\la)$ also differ at most by one, in addition to satisfy \eqref{eq.minfact2}.}
	
\section{Generic minimal rank factorizations in $\CC[\la]^{m\times n}_{d,r}$ and related results} \label{sect.compact}
In this section, {\blue for $r < \min \{m,n \}$,} we prove that arbitrarily close (in the distance defined in \eqref{ed.defdistance}) to any polynomial matrix $P(\la)\in \CC[\la]^{m\times n}_{d,r}$ there is another polynomial matrix $Q(\la)\in \CC[\la]^{m\times n}_{d,r}$ that can be factorized as $Q(\la) = L(\la) R(\la)$, with $L(\la) \in \CC[\la]^{m\times r}$, $R(\la) \in \CC[\la]^{r\times n}$, with the degrees of the columns of $L(\la)$ differing at most by one, with the degrees of the rows of $R(\la)$ also differing at most by one, {\blue  and satisfying \eqref{eq.minfact2}. Moreover, we prove that generically only $rd+1$ different degree distributions with such properties are necessary or, equivalently, we prove that $\CC[\la]^{m\times n}_{d,r}$ is the union of the closures of $rd+1$ sets of polynomial matrices which have rank-revealing factorizations with these very specific degree properties. In addition, we will see that $L(\la)$ and $R(\la)$ can be chosen to be minimal bases.} Finally, we will relate the sets of polynomial matrices in $\CC[\la]^{m\times n}_{d,r}$ that can be factorized in these specific ways with the orbits $\orb( K_a)$ of Theorem \ref{thm.mainthandrii2017}.

{\blue In contrast with Section \ref{sec.minfactorizations}, we assume throughout this section that $r<\min \{m, n\}$, i.e., the full rank case is not considered. The reason is that this section deals with generic results and, as explained in the last paragraph of Subsection \ref{subsec.genericandrii}, if $r= \min \{m, n\}$ and $m\ne n$, then generically the polynomial matrices in $\CC[\la]^{m\times n}_{d,r}$ are minimal bases with the degrees of their columns all equal to $d$ when $m>n$ or with the degrees of their rows all equal to $d$ when $m<n$. Therefore, they satisfy automatically \eqref{eq.minfact1} with one of the factors equal to the identity matrix and it makes no sense to look for minimal rank factorizations of these matrices. If $r= \min \{m, n\}$ and $m = n$, then generically the polynomial matrices in $\CC[\la]^{m\times n}_{d,r} = \CC[\la]^{n\times n}_{d}$ are regular with degree exactly $d$ (see again Subsection \ref{subsec.genericandrii}) and their minimal rank factorizations reduce to a triviality according to Remark \ref{rem.fullrank}.}

Before stating the first result in this section, we recall that the degree of the zero polynomial has been defined to be $-\infty$. Therefore, an expression as $\deg (L_{*i}) + \deg (R_{i*}) = d$ for the $i$th column and row of the factors in  $Q(\la) = L(\la) R(\la)$ implies that $L_{*i} (\la) \ne 0$, $R_{i*} (\la) \ne 0$, $0\leq \deg (L_{*i}) \leq d$, and $0\leq \deg (R_{i*}) \leq d$. In contrast, an expression as $\deg (L_{*i}) + \deg (R_{i*}) \leq d$ without further conditions does not imply that $\deg (L_{*i}) \leq d$ and $\deg (R_{i*}) \leq d$, because it might be possible that $\deg (L_{*i}) = -\infty$ and $\deg (R_{i*})$ is arbitrarily large, or vice versa.

The first result in this section is a simple consequence of the results in Section \ref{sec.minfactorizations} and states that every polynomial matrix in $\FF[\la]^{m\times n}_{d,r}$ can be factorized into two factors that reveal the maximum possible rank $r$ and such that the sums of the degrees of their corresponding columns and rows is bounded by $d$.

\begin{theo} \label{thm.1description} Let $m,n,r$ and $d$ be integers such that $m,n \geq 2$, $d \geq 1$ and $0 < r < \min\{m,n\}$. Then
\begin{equation} \label{eq.calS}
\FF[\la]^{m\times n}_{d,r} = \left\{ L(\la) R(\la) \, : \,
\begin{array}{l}
L(\la) \in \FF[\la]^{m\times r}, \; R(\la) \in \FF[\la]^{r\times n},\\
\deg (L_{*i}) \leq d, \quad \deg (R_{i*}) \leq d,\\
\deg (L_{*i}) + \deg (R_{i*}) \leq d, \quad \mbox{for $i=1,\ldots ,r$}
\end{array}
\right\} \, .
\end{equation}
\end{theo}

\begin{proof} {\blue Let $\mathcal{S}$ be the set in the right-hand side of \eqref{eq.calS}.}

Proof of $\FF[\la]^{m\times n}_{d,r} \subseteq \mathcal{S}$.
If $P(\la) \in \FF[\la]^{m\times n}_{d,r}$ and $P(\la) = 0$, then trivially $P(\la) = 0_{m\times r} 0_{r\times n} \in \mathcal{S}$. If $P(\la) \in \FF[\la]^{m\times n}_{d,r}$ and $P(\la) \ne 0$, then $0\leq \deg(P) = \widetilde{d} \leq d$ and $0 < \rank(P) = \widetilde{r} \leq r$. Then, Theorem \ref{thm.minfact0}-(ii) and {\blue Corollary \ref{thm.degreepredic}-(i)} imply that $P(\la)$ can be factorized as $P(\la) = \widetilde{L} (\la) \widetilde{R} (\la)$, with $\widetilde{L}(\la) \in \FF[\la]^{m\times \widetilde{r}}, \; \widetilde{R}(\la) \in \FF[\la]^{\widetilde{r}\times n}$, and
$0\leq \deg(\widetilde{L}_{*i}) + \deg(\widetilde{R}_{i*}) \leq \widetilde{d} \leq d$ for $i=1,\ldots , \widetilde{r}$. If $\widetilde{r} = r$, this proves that $P(\la) \in \mathcal{S}$. If $\widetilde{r} < r$, then we pad $\widetilde{L} (\la)$ and  $\widetilde{R} (\la)$ with zeros and define
$$
L(\la) := \left[\begin{array}{cc}
           \widetilde{L}(\la) & 0
         \end{array}\right] \in \FF[\la]^{m\times r} \quad \mbox{and} \quad
R(\la) := \left[\begin{array}{c}
           \widetilde{R}(\lambda) \\ 0
          \end{array}\right]  \in \FF[\la]^{r\times n}    \, ,
$$
which satisfy $P(\la) = L (\la) R (\la)$, with $\deg (L_{*i}) \leq d$, $\deg (R_{i*}) \leq d$, and $\deg(L_{*i}) + \deg(R_{i*}) \leq d$ for $i=1,\ldots , r$. Therefore, $P(\la) \in \mathcal{S}$. This proves $\FF[\la]^{m\times n}_{d,r} \subseteq \mathcal{S}$.

Proof of $\mathcal{S}  \subseteq \FF[\la]^{m\times n}_{d,r}$.
If $P(\la) = L(\la) R(\la) \in \mathcal{S}$, then $\rank(P) \leq \min \{ \rank(L) , \rank(R)\}\allowbreak \leq r$. In addition, the expansion $P(\la) = \sum_{i=1}^r L_{*i}(\la) R_{i*} (\la)$ and Lemma \ref{lemm.degreerank1} imply $\deg(P) \leq \max_{1 \leq i \leq r} \{ \deg(L_{*i}) + \deg ( R_{i*} ) \} \leq d$. Thus $P(\la) \in \FF[\la]^{m\times n}_{d,r}$, and the proof is completed.
\end{proof}

The rest of the results of this section are valid only over the field $\CC$ since they use limits and topological concepts with respect to the distance in \eqref{ed.defdistance}. This will allow us to prove that every polynomial matrix in $\CC[\la]^{m\times n}_{d,r}$ is the limit of a sequence of polynomial matrices in $\CC[\la]^{m\times n}_{d,r}$ that can be factorized into two factors such that the degrees of their columns and rows  have very specific properties when compared with those in Theorem \ref{thm.1description}. The first result in this direction is Theorem \ref{thm.2description}.

\begin{theo} \label{thm.2description} Let $m,n,r$ and $d$ be integers such that $m,n \geq 2$, $d \geq 1$ and $0 < r < \min\{m,n\}$ and define the sets
$$
\mathcal{A}^{m\times n}_{d,r} := \left\{ L(\la) R(\la) \, : \,
\begin{array}{l}
L(\la) \in \CC[\la]^{m\times r}, \; R(\la) \in \CC[\la]^{r\times n},\\
\deg (L_{*i}) + \deg (R_{i*}) = d, \quad \mbox{for $i=1,\ldots ,r$}
\end{array}
\right\} .
$$
Then
$$
\CC[\la]^{m\times n}_{d,r} = \overline{\mathcal{A}^{m\times n}_{d,r}}  \, .
$$
\end{theo}

\begin{proof}
From Theorem \ref{thm.1description} it is obvious that $\mathcal{A}^{m\times n}_{d,r} \subseteq \CC[\la]^{m\times n}_{d,r}$. Moreover, $\CC[\la]^{m\times n}_{d,r}$ is a closed subset of $\CC[\la]^{m\times n}_d$ and the closure of $\mathcal{A}^{m\times n}_{d,r}$ is the smallest closed set that contains $\mathcal{A}^{m\times n}_{d,r}$. Therefore, $\mathcal{A}^{m\times n}_{d,r} \subseteq \overline{\mathcal{A}^{m\times n}_{d,r}} \subseteq \CC[\la]^{m\times n}_{d,r}$.

In the rest of the proof, we prove that $\CC[\la]^{m\times n}_{d,r} \subseteq \overline{\mathcal{A}^{m\times n}_{d,r}}$. If $P(\la) \in \CC[\la]^{m\times n}_{d,r}$, then Theorem \ref{thm.1description} implies that $P(\la) = L(\la) R(\la)$ with $L(\la) \in \CC[\la]^{m\times r}, \; R(\la) \in \CC[\la]^{r\times n}$, $\deg (L_{*i}) \leq d$, $\deg (R_{i*}) \leq d$, and
$\deg (L_{*i}) + \deg (R_{i*}) \leq d,$ for $i=1,\ldots ,r$. Moreover, without loss of generality, we take $L_{*j} (\la) = 0$ whenever $R_{
j*} (\la) = 0$. If $\deg (L_{*i}) + \deg (R_{i*}) = d,$ for $i=1,\ldots ,r$, then  $P(\la) \in \mathcal{A}^{m\times n}_{d,r}  \subseteq \overline{\mathcal{A}^{m\times n}_{d,r} }$. Otherwise, let us consider the set of indices corresponding to strict inequalities, that is,
$$\mathcal{I} := \{ j \, : \, 1 \leq j \leq r \; \mbox{and} \; \deg (L_{*j}) + \deg (R_{j*}) < d\} .$$
Then, consider any two sequences of constant nonzero vectors $\{v_k\}_{k\in \NN} \subset \CC^{m\times 1}$ such that $\lim_{k\rightarrow \infty} v_k = 0$ and $\{w_k\}_{k\in \NN} \subset \CC^{1\times n}$ such that $\lim_{k\rightarrow \infty} w_k = 0$, and construct the following two sequences of polynomial matrices: (1) $L_k (\la) = L (\la) + F_k (\la)$, where the columns of $F_k (\la)$ are constructed as follows
$$(F_k)_{*j} (\la) = \left\{ \begin{array}{ll}
                            0, &  \mbox{if $j \notin \mathcal{I}$}, \\
                            \la^{d - \deg (R_{j*})} v_k, & \mbox{if $j \in \mathcal{I}$ and $R_{j*}(\la) \ne 0$}, \\
                            \la^{d} v_k, & \mbox{if $j \in \mathcal{I}$ and $R_{j*}(\la) = 0$},
                          \end{array} \right.$$
and (2) $R_k (\la) = R (\la) + G_k (\la)$, where the rows of $G_k (\la)$ are constructed as follows
$$(G_k)_{j*} (\la) = \left\{ \begin{array}{ll}
                            0, & \mbox{if $j \notin \mathcal{I}$}, \\
                            0, & \mbox{if $j \in \mathcal{I}$ and $R_{j*}(\la) \ne 0$}, \\
                            w_k, & \mbox{if $j \in \mathcal{I}$ and $R_{j*}(\la) = 0$}.
                          \end{array} \right.$$
Then, $P_k (\la) := L_k (\la) R_k(\la) \in \mathcal{A}^{m\times n}_{d,r}$ and $\lim_{k\rightarrow \infty}P_k(\la) = P(\la)$, which implies that $P(\la) \in \overline{\mathcal{A}^{m\times n}_{d,r}}$.
\end{proof}

Next, we consider some subsets of the set $\mathcal{A}^{m\times n}_{d,r}$ introduced in Theorem \ref{thm.2description} that will be fundamental auxiliary tools for getting the main results of this section. More precisely, we express in the next theorem the set $\mathcal{A}^{m\times n}_{d,r}$ as the union of such subsets, and $\CC[\la]^{m\times n}_{d,r}$ as the union of their closures.

\begin{theo} \label{thm.3description} Let $\mathcal{A}^{m\times n}_{d,r}$ be the set defined in Theorem \ref{thm.2description} and for each natural number $a = 0,1,....,rd$ define the following subsets of $\CC[\la]^{m\times n}_{d,r}$
$$ {\blue
\mathcal{A}^{m\times n}_{d,r,a} :=  \left\{ L(\la) R(\la) \, : \,
\begin{array}{l}
L(\la) R(\la) \in \mathcal{A}^{m\times n}_{d,r},\\
\sum_{i=1}^{r} \deg (R_{i*}) = a
\end{array}
\right\}  \, .}
$$
Then
\begin{itemize}
  \item[\rm (i)] $\displaystyle \mathcal{A}^{m\times n}_{d,r} = \bigcup_{0\leq a \leq rd} \mathcal{A}^{m\times n}_{d,r,a} $,
  \item[\rm (ii)] $\displaystyle \CC[\la]^{m\times n}_{d,r} = \bigcup_{0\leq a \leq rd} \overline{\mathcal{A}^{m\times n}_{d,r,a}}$,

  \item[\rm (iii)] for every $P(\la) \in \CC[\la]^{m\times n}_{d,r}$, there exists an integer $a$ such that $P(\la) \in \overline{\mathcal{A}^{m\times n}_{d,r,a}}$.
\end{itemize}
\end{theo}

\begin{proof}
Item (i). Let us prove first that $\mathcal{A}^{m\times n}_{d,r} \subseteq \bigcup_{0\leq a \leq rd} \mathcal{A}^{m\times n}_{d,r,a}$. If $P(\la) \in \mathcal{A}^{m\times n}_{d,r}$, then $P(\la) = L(\la) R(\la)$ with $L(\la) \in \CC[\la]^{m\times r}, \; R(\la) \in \CC[\la]^{r\times n}$, and
$\deg (L_{*i}) + \deg (R_{i*}) = d,$ for $i=1,\ldots ,r$. So, $\sum_{i=1}^{r}\deg (L_{*i}) + \sum_{i=1}^{r}\deg (R_{i*}) = rd$, which implies that
$0 \leq \sum_{i=1}^{r}\deg (R_{i*}) \leq rd$. Therefore, $P(\la) \in \mathcal{A}^{m\times n}_{d,r,a}$ for some $a =0,1,\ldots , rd$ and $P(\la) \in \bigcup_{0\leq a \leq rd} \mathcal{A}^{m\times n}_{d,r,a}$.

The reverse inclusion $\bigcup_{0\leq a \leq rd} \mathcal{A}^{m\times n}_{d,r,a} \subseteq \mathcal{A}^{m\times n}_{d,r}$ holds by definition.

\medskip

Item (ii). It is an immediate consequence of Theorem \ref{thm.2description}, item (i), and the basic fact that ``the closure of the union of a finite number of sets is the union of the closures of such sets''.

Item (iii) is just another expression of item (ii).
\end{proof}

We present next the key technical result of this section, Theorem \ref{thm.subsets}, which deals with some subsets of $\mathcal{A}^{m\times n}_{d,r,a}$ that are introduced in Definition \ref{def.rowdegsubsets}.

\begin{deff} \label{def.rowdegsubsets}
Let $m,n,r$ and $d$ be integers such that $m,n \geq 2$, $d \geq 1$ and $0 < r < \min\{m,n\}$, $a$ be an integer such that $0 \leq a \leq rd$, {\blue $(\rho_1, \ldots, \rho_r)$ be any list of integers such that $0\leq \rho_i \leq d$, for $i=1, 2 ,\ldots , r$, and $\sum_{i=1}^{r} \rho_i = a$, and $\mathcal{A}^{m\times n}_{d,r,a}$ be the set defined in Theorem \ref{thm.3description}.  The following subsets of polynomial matrices are defined
$$
\mathcal{A}^{m\times n}_{d,r,a} (\rho_1, \ldots, \rho_r) :=  \left\{ L(\la) R(\la)  :
\begin{array}{l}
L(\la) R(\la) \in \mathcal{A}^{m\times n}_{d,r,a},\\
\deg (R_{i*}) = \rho_i, \quad \deg (L_{*i}) = d - \rho_i, \; \mbox{for $i=1,\ldots ,r$}
\end{array} \! \!
\right\} .
$$}
\end{deff}



\begin{theo} \label{thm.subsets} Let $\mathcal{A}^{m\times n}_{d,r,a}$ be the set defined in Theorem \ref{thm.3description} and $\mathcal{A}^{m\times n}_{d,r,a} (\rho_1, \ldots, \rho_r)$ be any of the sets defined in Definition \ref{def.rowdegsubsets}. Then the following statements hold:
\begin{itemize}
  \item[\rm (i)] If $(\sigma_1, \ldots, \sigma_r)$ is any permutation of $(1,\ldots ,r)$, then
      $$
      \mathcal{A}^{m\times n}_{d,r,a} (\rho_1, \ldots, \rho_r) = \mathcal{A}^{m\times n}_{d,r,a} (\rho_{\sigma_1}, \ldots, \rho_{\sigma_r}).
      $$
  \item[\rm (ii)] If $\rho_j - \rho_k \geq 2$, then
      $$
      \mathcal{A}^{m\times n}_{d,r,a} (\rho_1, \ldots, \rho_j , \ldots , \rho_k,\ldots, \rho_r)
      \subseteq
      \overline{\mathcal{A}^{m\times n}_{d,r,a} (\rho_1, \ldots, \rho_j -1  , \ldots , \rho_k +1 ,\ldots, \rho_r)} \, .
      $$

  \item[\rm (iii)] If $d_R = \left\lfloor a /r \right\rfloor$ and $t_R = a \; \mathrm{mod} \; r$, then
      $$
      \mathcal{A}^{m\times n}_{d,r,a} (\rho_1, \ldots, \rho_r) \subseteq
      \overline{\mathcal{A}^{m\times n}_{d,r,a} (\underbrace{d_R +1, \ldots, d_R + 1}_{t_R}  , \underbrace{d_R,\ldots, d_R}_{r-t_R})} \, .
      $$
\end{itemize}

\end{theo}

\begin{proof}
Proof of item (i). If $L(\la) R(\la) \in \mathcal{A}^{m\times n}_{d,r,a} (\rho_1, \ldots, \rho_r)$ and $\Pi$ is an $r\times r$ permutation matrix such that the $i$th row of $\Pi R(\la)$ is the $\sigma_i$th row of $R(\la)$, for $i=1,\ldots ,r$, then $L(\la) R(\la) = (L(\la) \Pi^T) (\Pi R(\la) ) \in \mathcal{A}^{m\times n}_{d,r,a} (\rho_{\sigma_1}, \ldots, \rho_{\sigma_r})$. Therefore, $\mathcal{A}^{m\times n}_{d,r,a} (\rho_1, \ldots, \rho_r) \subseteq \mathcal{A}^{m\times n}_{d,r,a} (\rho_{\sigma_1}, \ldots, \rho_{\sigma_r})$. The ``reverse'' inclusion is proved in a similar manner using the ``reverse'' permutation.

Proof of item (ii). As a consequence of item (i), we can assume without loss of generality that $j=1$ and $k=2$. Let $P(\la) = L(\la) R(\la) \in \mathcal{A}^{m\times n}_{d,r,a} (\rho_1 , \rho_2 , \ldots, \rho_r)$ with $\rho_1 - \rho_2 \geq 2$. Then, the first row of $R(\la)$ and the second column of $L(\la)$ can be written as follows:
\begin{align} \label{eq.auxepsilonk1}
R_{1*} (\la) & = \la^{\rho_1} v_{\rho_1} +  \widetilde{R}_{1*} (\la), \phantom{iiiiiiiii}\mbox{with $0 \ne v_{\rho_1}\in \CC^{1 \times n}$ and $\deg(\widetilde{R}_{1*}) < \rho_1$},\\ \label{eq.auxepsilonk2}
L_{*2} (\la) & = \la^{d-\rho_2} w_{d-\rho_2} +  \widetilde{L}_{*2} (\la), \quad \mbox{with $0 \ne w_{d-\rho_2}\in \CC^{m\times 1}$ and $\deg(\widetilde{L}_{*2}) < d-\rho_2$}.
\end{align}
Next, for any sequence $\{\epsilon_k \}_{k\in \NN} \subset \CC$ of {\it nonzero} numbers such that $\lim_{k \rightarrow \infty} \epsilon_k = 0$, we define two sequences of polynomial matrices $\{L_k (\la)\}_{k\in \NN} \subseteq \CC[\la]^{m\times r}$ and
$\{R_k (\la)\}_{k\in \NN} \subseteq \CC[\la]^{r\times n}$ (via their columns and rows, respectively) as follows
\begin{equation} \label{eq.auxepsilonk3} \!\!\!
\begin{array}{l}
(L_k)_{*1} (\la) := -\epsilon_k \la^{d - \rho_1 + 1} w_{d-\rho_2} + L_{*1} (\la), \\
(L_k)_{*i} (\la) := L_{*i} (\la), \quad \mbox{$1< i \leq r$},
\end{array}
\begin{array}{l}
(R_k)_{2*} (\la) := \epsilon_k \la^{\rho_2 + 1} v_{\rho_1} + R_{2*} (\la), \\
(R_k)_{i*} (\la) := R_{i*} (\la), \quad \mbox{$i\ne 2$, $1\leq i \leq r$}.
\end{array}
\end{equation}
From these sequences, we define the sequence $\{P_k(\la)\}_{k\in \NN} := \{ L_k (\la) R_k (\la)\}_{k\in \NN} \subseteq \CC[\la]^{m\times n}$, which obviously satisfies $\lim_{k\rightarrow \infty} P_k(\la) = P(\la)$. In the rest of the proof, we will prove that there exists an index $k_0$ such that for every $k \geq k_0$, $$P_k(\la) = L_k (\la) R_k (\la) \in \mathcal{A}^{m\times n}_{d,r,a} (\rho_1 -1, \rho_2 + 1,\rho_3,\ldots, \rho_r),$$ which implies that $P(\la) \in
\overline{\mathcal{A}^{m\times n}_{d,r,a} (\rho_1 -1, \rho_2 + 1,\rho_3, \ldots, \rho_r)}$.
For this purpose, we define
$$D_k (\la) := \left[\begin{array}{cc}
            1 & -\frac{1}{\epsilon_k} \, \la^{\rho_1 - \rho_2 -1 } \\
            0 & 1
          \end{array} \right] \oplus I_{r-2},$$
whose inverse is
$$
D_k (\la)^{-1} := \left[\begin{array}{cc}
            1 & \frac{1}{\epsilon_k} \, \la^{\rho_1 - \rho_2 -1 } \\
            0 & 1
          \end{array} \right] \oplus I_{r-2} \, .$$
Therefore,
\begin{equation} \label{eq.auxclos1}
P_k (\la) = (L_k (\la) D_k (\la)^{-1}) (D_k (\la) R_k (\la)).
\end{equation}
The $i$th row of $(D_k (\la) R_k (\la))$ is equal to the $i$th row of $R_k (\la)$ for $i = 2, \ldots, r$, and, taking into account \eqref{eq.auxepsilonk1} and \eqref{eq.auxepsilonk3}, the first row is
$$
(D_k (\la) R_k (\la))_{1*} = \widetilde{R}_{1*} (\la) - \frac{1}{\epsilon_k} \la^{\rho_1 - \rho_2 -1 } R_{2*} (\la),
$$
which has degree $\rho_1 -1$ for $\epsilon_k$ sufficiently close to zero or equivalently for all $k$ sufficiently large. In summary, there exists an index $k'$ such that for all $k\geq k'$
\begin{equation} \label{eq.auxclos2}
\mbox{the degrees of the rows of $D_k (\la) R_k (\la)$ are} \quad
\rho_1-1, \rho_2+1, \rho_3 , \rho_4 , \ldots , \rho_r \, .
\end{equation}
On the other hand, the $i$th column of $L_k (\la) D_k (\la)^{-1}$ is equal to the $i$th column of $L_k (\la)$ for $i=1,3,4 \ldots , r$, and, taking into account \eqref{eq.auxepsilonk2} and \eqref{eq.auxepsilonk3}, the second column is
$$
(L_k (\la) D_k (\la)^{-1})_{*2} = \frac{1}{\epsilon_k} \la^{\rho_1 - \rho_2 -1} \, L_{*1} (\la) + \widetilde{L}_{*2} (\la) ,
$$
which has degree $d - \rho_2 -1$ for $\epsilon_k$ sufficiently close to zero or equivalently for all $k$ sufficiently large. In summary, there exists an index $k''$ such that for all $k\geq k''$
\begin{equation} \label{eq.auxclos3}
\mbox{the degrees of the columns of $L_k (\la) D_k (\la)^{-1}$ are} \quad
d - \rho_1 +1, d - \rho_2 -1, d - \rho_3, d-\rho_4, \ldots , d - \rho_r \,.
\end{equation}
Combining \eqref{eq.auxclos1}, \eqref{eq.auxclos2}, and \eqref{eq.auxclos3} we get that $$P_k(\la) = L_k (\la) R_k (\la) \in \mathcal{A}^{m\times n}_{d,r,a} (\rho_1 -1, \rho_2 + 1,\rho_3, \ldots, \rho_r)$$ for all $k\geq \max \{k', k''\} = k_0$ and the proof is completed.

Proof of item (iii). Observe that item (ii) and the fact that ``the closure of a set is the smallest closed set that includes it'' imply
$$
\overline{\mathcal{A}^{m\times n}_{d,r,a} (\rho_1, \ldots, \rho_j , \ldots , \rho_k,\ldots, \rho_r)}
\subseteq
\overline{\mathcal{A}^{m\times n}_{d,r,a} (\rho_1, \ldots, \rho_j -1  , \ldots , \rho_k +1 ,\ldots, \rho_r)} \, .
$$
Therefore, we can apply again this result to the set on the right hand side of the equation above (permuting if necessary the indices by using the result in item (i)) as long as for at least two of the indices in $(\rho_1, \ldots, \rho_j -1  , \ldots , \rho_k +1 ,\ldots, \rho_r)$ the absolute value of their difference is larger than or equal to two. We can construct in this way a chain of subset inclusions until the indices $\rho_i$ differ at most by one unit, that is,
\begin{align*}
\mathcal{A}^{m\times n}_{d,r,a} (\rho_1, \ldots, \rho_j , \ldots , \rho_k,\ldots, \rho_r) &
\subseteq
\overline{\mathcal{A}^{m\times n}_{d,r,a} (\rho_1, \ldots, \rho_j , \ldots , \rho_k,\ldots, \rho_r)} \\ &
\subseteq
\overline{\mathcal{A}^{m\times n}_{d,r,a} (\rho_1, \ldots, \rho_j -1  , \ldots , \rho_k +1 ,\ldots, \rho_r)} \\
& \subseteq \cdots \subseteq  \overline{\mathcal{A}^{m\times n}_{d,r,a} (\underbrace{d_R +1, \ldots, d_R + 1}_{t_R}  , \underbrace{d_R,\ldots, d_R}_{r-t_R})} \, .
\end{align*}
We emphasize that the values of $d_R$ and $t_R$ are completely determined by the fact that the sum of the $r$ indices of all the subsets in the chain above is always $a$ and that the indices in the last subset differ at most by one (in absolute value).
\end{proof}

\begin{example} \label{ex.technicalthm} {\rm In order to illustrate the proof and the statement of Theorem \ref{thm.subsets}, we consider the following polynomial matrix
\begin{equation}\label{eq.extechnicalthm}
P (\la) =
\left[
\begin{array}{cc}
  0 & \la^2 \\
  1 & 1 \\
  1 & 0
\end{array}
\right]
\left[
\begin{array}{ccc}
  0 & \la^2 & 1\\
  1 & 0 & 0
\end{array}
\right] =
\left[
\begin{array}{ccc}
\la^2 & 0 & 0 \\
 1 & \la^2 & 1\\
 0 & \la^2 & 1
\end{array}
\right] \in \mathcal{A}^{3\times 3}_{2,2,2} (2,0) \subset \CC[\la]_{2,2}^{3 \times 3} \, .
\end{equation}
Since $a=2$ and $r=2$, the quantities in Theorem \ref{thm.subsets}-(iii) are $d_R =1$ and $t_R =0$. One might wonder whether $P(\la)$ might be factorized in a form different from the one in \eqref{eq.extechnicalthm} in such a way that $P(\la)\in \mathcal{A}^{3\times 3}_{2,2,2} (1,1)$. However, it is easy to see that $P(\la)\notin \mathcal{A}^{3\times 3}_{2,2,2} (1,1)$ as follows. Observe first that in the factorization $P(\la) = L(\la) R(\la)$ given in \eqref{eq.extechnicalthm} both factors are minimal bases by Theorem \ref{thm.minbasischar}. Thus, the minimal indices of ${\mathcal Row}(P)$ are $2$ and $0$. If $P(\la)\in \mathcal{A}^{3\times 3}_{2,2,2} (1,1)$, then there would exist a factorization $P(\la) = \widetilde{L}(\la) \widetilde{R}(\la)$ with $\widetilde{L}(\la) \in \CC[\la]^{3\times 2}$ and $\widetilde{R}(\la) \in \CC[\la]^{2\times 3}$ with the degrees of both rows of $\widetilde{R}(\la)$ equal to $1$ and, since $\rank (P) = 2$, Lemma \ref{lemm:fact0}-(iv) would imply that the rows of $\widetilde{R}(\la)$ form a polynomial basis of ${\mathcal Row}(P)$ with the sum of the degrees of its vectors equal to $2$. Therefore, the rows of $\widetilde{R}(\la)$ would be a minimal basis of ${\mathcal Row}(P)$ and the minimal indices of this rational subspace would be $1$ and $1$, which contradicts that the minimal indices of ${\mathcal Row}(P)$ are $2$ and $0$.

Consider any sequence $\{\epsilon_k\}_{k\in \NN}$ of nonzero numbers with $\lim_{k\rightarrow \infty} \epsilon_k = 0$ and construct from $P(\la)$ the following sequence of polynomial matrices via the strategy in \eqref{eq.auxepsilonk3}:
$$
P_k (\la) =
\left[
\begin{array}{cc}
  -\epsilon_k \la & \la^2 \\
  1 & 1 \\
  1 & 0
\end{array}
\right]
\left[
\begin{array}{ccc}
  0 & \la^2 & 1\\
  1 & \epsilon_k \la & 0
\end{array}
\right] =
\left[
\begin{array}{ccc}
\la^2 & 0 & -\epsilon_k \la \\
 1 & \la^2 + \epsilon_k \la & 1\\
 0 & \la^2 & 1
\end{array}
\right] ,
$$
which satisfies $\lim_{k\rightarrow \infty} P_k (\la) = P(\la)$.
Proceeding as in \eqref{eq.auxclos1}, $P_k (\la)$ can be written as:
\begin{align*}
P_k (\la) & =
\left[
\begin{array}{cc}
  -\epsilon_k \la & \la^2 \\
  1 & 1 \\
  1 & 0
\end{array}
\right]
\left[
\begin{array}{cc}
  1 & \frac{1}{\epsilon_k} \, \lambda \\
  0 & 1
\end{array}
\right]
\left[
\begin{array}{cc}
  1 & -\frac{1}{\epsilon_k} \, \lambda \\
  0 & 1
\end{array}
\right]
\left[
\begin{array}{ccc}
  0 & \la^2 & 1\\
  1 & \epsilon_k \la & 0
\end{array}
\right] \\
& =
\left[
\begin{array}{cc}
  -\epsilon_k \la & 0 \\
  1 & \frac{1}{\epsilon_k} \la + 1 \\
  1 & \frac{1}{\epsilon_k} \la
\end{array}
\right]
\left[
\begin{array}{ccc}
  -\frac{1}{\epsilon_k} \la & 0 & 1\\
                        1 & \epsilon_k \la & 0
\end{array}
\right] \in \mathcal{A}^{3\times 3}_{2,2,2} (1,1) \subset \CC[\la]_{2,2}^{3 \times 3} \,.
\end{align*} {\blue \qed}
}
\end{example}

{\blue The set in Theorem \ref{thm.subsets}-(iii) plays a crucial role in the main results of this section. Therefore, we adopt the following short notation for it $\mathcal{B}^{m\times n}_{d,r,a} = \mathcal{A}^{m\times n}_{d,r,a} (d_R +1, \ldots, d_R + 1  , d_R,\ldots, d_R)$, and we define it explicitly in Definition \ref{def.homogsubsets} for future reference.}

\begin{deff} \label{def.homogsubsets} Let $m,n,r$ and $d$ be integers such that $m,n \geq 2$, $d \geq 1$ and $0 < r < \min\{m,n\}$, and $a$ be an integer such that $0 \leq a \leq rd$. Let us define $d_R := \left\lfloor a /r \right\rfloor$, $t_R := a \; \mathrm{mod} \; r$ and the following subset of polynomial matrices
$$
\mathcal{B}^{m\times n}_{d,r,a}  :=  \left\{ L(\la) R(\la) \, : \,
\begin{array}{l}
L(\la) \in \CC[\la]^{m\times r}, \; R(\la) \in \CC[\la]^{r\times n},\\
\deg (R_{i*}) = d_R + 1, \quad \mbox{for $i=1,\ldots ,t_R$}, \\
\deg (R_{i*}) = d_R, \quad \mbox{for $i=t_R + 1,\ldots ,r$}, \\
\deg (L_{*i}) = d - \deg (R_{i*}), \quad \mbox{for $i=1,\ldots ,r$}
\end{array}
\right\}  \subset \CC[\la]^{m\times n}_{d,r} \, .
$$
\end{deff}


As a simple consequence of the developments above, we prove the first main result of this section.
\begin{theo} \label{thm.4descrition} Let $\mathcal{A}^{m\times n}_{d,r,a}$ and $\mathcal{B}^{m\times n}_{d,r,a}$ be the sets of polynomial matrices defined in Theorem \ref{thm.3description} and in Definition \ref{def.homogsubsets}, respectively. Then,
\begin{itemize}
\item[\rm (i)] $\mathcal{B}^{m\times n}_{d,r,a} \subseteq \mathcal{A}^{m\times n}_{d,r,a}$  for $a = 0,1,\ldots , rd$,

\item[\rm (ii)] $\overline{\mathcal{B}^{m\times n}_{d,r,a}} = \overline{\mathcal{A}^{m\times n}_{d,r,a}}$ for $a = 0,1,\ldots , rd$,

\item[\rm (iii)] $\displaystyle \CC[\la]^{m\times n}_{d,r} = \bigcup_{0 \leq a \leq rd} \overline{\mathcal{B}^{m\times n}_{d,r,a}} \, ,$ and

\item[\rm (iv)] for every $P(\la) \in \CC[\la]^{m\times n}_{d,r}$, there exists an integer $a$ such that $P(\la) \in \overline{\mathcal{B}^{m\times n}_{d,r,a}}$.

\end{itemize}
\end{theo}

\begin{proof} Item (i) is obvious by definition. Item (i) implies  $\overline{\mathcal{B}^{m\times n}_{d,r,a}} \subseteq \overline{\mathcal{A}^{m\times n}_{d,r,a}}$. Next, suppose $L(\la) R(\la) \in \mathcal{A}^{m\times n}_{d,r,a}$. Then $L(\la) R(\la) \in \mathcal{A}^{m\times n}_{d,r,a} (\rho_1,\ldots, \rho_r)$ for some integers $(\rho_1,\ldots, \rho_r)$ such that $0\leq \rho_i \leq d$, for $i= 1, \ldots ,r,$ and $\sum_{i=1}^{r} \rho_i = a$, and, by Theorem \ref{thm.subsets}-(iii), $L(\la) R(\la) \in \overline{\mathcal{B}^{m\times n}_{d,r,a}}$. Therefore, $\mathcal{A}^{m\times n}_{d,r,a} \subseteq \overline{\mathcal{B}^{m\times n}_{d,r,a}}$, which implies $\overline{\mathcal{A}^{m\times n}_{d,r,a}} \subseteq \overline{\mathcal{B}^{m\times n}_{d,r,a}}$. This proves item (ii). Finally, items (iii) and (iv) follow from item (ii) and the items (ii) and (iii), respectively, of Theorem \ref{thm.3description}.
\end{proof}

{\blue We know that there are polynomial matrices for which none of their rank factorizations satifies that the sum of the degrees of the factors is equal to the degree of the polynomial. Recall Lemma \ref{lemm:lowbounddeg} and Example \ref{ex.arbitrarydeg}. However, a corollary of Theorem \ref{thm.4descrition} is that generically the polynomials in $\CC[\la]^{m\times n}_{d,r}$ have factorizations with the sum of the degrees of the factors not larger than $d+1$. The reason is that if $L(\la) R(\la) \in \mathcal{B}^{m\times n}_{d,r,a}$, then $\deg (L) + \deg (R) = d +1$ if $t_R >0$ and $\deg (L) + \deg (R) = d$ if $t_R = 0$. We state this result as a corollary for future reference.
\begin{cor} Let $${\mathcal S}_{d,r}^{m \times n}:=
\left\{ L(\la) R(\la) \, : \,
\begin{array}{l}
L(\la) \in \CC[\la]^{m\times r}, \; R(\la) \in \CC[\la]^{r\times n},\\
\deg (L) + \deg (R) \leq d + 1, \\
L(\la) R(\la) \in \CC[\la]^{m \times n}_{d,r}
\end{array}
\right\}
\subset \CC[\la]^{m \times n}_{d,r}.$$ Then $\CC[\la]^{m \times n}_{d,r} = \overline{{\mathcal S}_{d,r}^{m \times n}}$.
\end{cor}

}

Theorem \ref{thm.4descrition} proves the promised result that arbitrarily close to any polynomial matrix $P(\la)\in \CC[\la]^{m\times n}_{d,r}$, there is another polynomial matrix $Q(\la)\in \CC[\la]^{m\times n}_{d,r}$ that can be factorized as $Q(\la) = L(\la) R(\la)$, with the degrees of the columns of $L(\la)$ differing at most by one, with the degrees of the rows of $R(\la)$ also differing at most by one, {\blue and, for each $i=1, \ldots, r$, the sum of the degrees of the $i$th column of $L(\la)$ and of the $i$th row of $R(\lambda)$ is equal to $d$. Moreover, we have proved that only $rd+1$ different degree distributions with these properties are necessary}. However, the factorization of $Q(\la)$ is not necessarily a minimal rank factorization, according to the definition of $\mathcal{B}^{m\times n}_{d,r,a}$. Next, we prove in Theorem \ref{thm.5descrition} that arbitrarily close to any polynomial matrix $P(\la)\in \CC[\la]^{m\times n}_{d,r}$ there is a polynomial matrix $Q(\la)$ that can be factorized as $Q(\la)= L(\la) R(\la)$ with factors  satisfying the conditions of Theorem \ref{thm.minfact}, and, moreover, with the degrees of the columns of $L(\la)$ differing at most by one and with the degrees of the rows of $R(\la)$ also differing at most by one. In addition, the minimal indices of ${\mathcal N}_\ell (Q)$ and ${\mathcal N}_r (Q)$ are as those in Theorem \ref{thm.mainthandrii2017}. For that purpose we introduce first the following definitions.

\begin{deff} \label{def.minbassubsets} Let $m,n,r$ and $d$ be integers such that $m,n \geq 2$, $d \geq 1$ and $0 < r < \min\{m,n\}$, {\blue $a$ be an integer such that $0 \leq a \leq rd$, and $\mathcal{B}^{m\times n}_{d,r,a}$ be the set in Definition \ref{def.homogsubsets}. Let us define $\alpha := \lfloor a / (n-r) \rfloor$, $s := a \mod (n-r)$, $\beta := \lfloor (rd-a)/(m-r) \rfloor$,
and $t := (rd-a) \mod (m-r)$ and the following subsets of $\CC[\la]^{m\times n}_{d,r}$} {\blue
$$
\mathcal{M}^{m\times n}_{d,r,a}  :=  \left\{ L(\la) R(\la) \, : \,
\begin{array}{l}
L(\la)  R(\la) \in \mathcal{B}^{m\times n}_{d,r,a},\\
L(\la) \; \mbox{and} \; R(\la) \; \mbox{are minimal bases}
\end{array}
\right\}  \, ,
$$
$$
\mathcal{MH}^{m\times n}_{d,r,a}  :=  \left\{ L(\la) R(\la) \, : \,
\begin{array}{l}
L(\la) R(\la) \in \mathcal{M}^{m\times n}_{d,r,a},\\
\mathcal{N}_\ell (L) \; \mbox{has minimal indices} \;
\{\underbrace{\beta+1, \dots , \beta+1}_{t}, \underbrace{\beta, \dots , \beta}_{m-r-t}\} ,\\
\mathcal{N}_r (R) \; \mbox{has minimal indices} \;
\{\underbrace{\alpha+1, \dots , \alpha+1}_{s},\underbrace{\alpha, \dots , \alpha}_{n-r-s}\}
\end{array}
\right\}  \, .
$$}
\end{deff}

With respect to the definition of $\mathcal{MH}^{m\times n}_{d,r,a}$, it is important to recall that Lemma \ref{lemm:fact0} implies that $\mathcal{N}_\ell (L) = \mathcal{N}_\ell (L(\la) R(\la))$ and that $\mathcal{N}_r (R) = \mathcal{N}_r (L(\la) R(\la))$.

\begin{theo} \label{thm.5descrition} Let $\mathcal{B}^{m\times n}_{d,r,a}$, $\mathcal{M}^{m\times n}_{d,r,a}$ and $\mathcal{MH}^{m\times n}_{d,r,a}$ be the sets of polynomial matrices introduced in Definitions \ref{def.homogsubsets} and \ref{def.minbassubsets}. Then,
\begin{itemize}
\item[\rm (i)] $\mathcal{MH}^{m\times n}_{d,r,a} \subseteq \mathcal{M}^{m\times n}_{d,r,a} \subseteq \mathcal{B}^{m\times n}_{d,r,a}$  for $a = 0,1,\ldots , rd$,

\item[\rm (ii)] $\overline{\mathcal{MH}^{m\times n}_{d,r,a}}=\overline{\mathcal{M}^{m\times n}_{d,r,a}} = \overline{\mathcal{B}^{m\times n}_{d,r,a}}$ for $a = 0,1,\ldots , rd$,

\item[\rm (iii)] $\displaystyle \CC[\la]^{m\times n}_{d,r} = \bigcup_{0 \leq a \leq rd} \overline{\mathcal{MH}^{m\times n}_{d,r,a}} = \bigcup_{0 \leq a \leq rd} \overline{\mathcal{M}^{m\times n}_{d,r,a}} \, ,$ and

\item[\rm (iv)] for every $P(\la) \in \CC[\la]^{m\times n}_{d,r}$, there exists an integer $a$ such that $P(\la) \in \overline{\mathcal{MH}^{m\times n}_{d,r,a}} = \overline{\mathcal{M}^{m\times n}_{d,r,a}}$.

\end{itemize}
\end{theo}

\begin{proof}
Item (i) is obvious from the definitions of the involved sets.

Proof of item (ii). First note that item (i) implies immediately that $$\overline{\mathcal{MH}^{m\times n}_{d,r,a}} \subseteq \overline{\mathcal{M}^{m\times n}_{d,r,a}} \subseteq \overline{\mathcal{B}^{m\times n}_{d,r,a}}.$$ With this result at hand, observe that if we prove $\mathcal{B}^{m\times n}_{d,r,a} \subseteq \overline{\mathcal{MH}^{m\times n}_{d,r,a}}$, then $\overline{\mathcal{B}^{m\times n}_{d,r,a}} \subseteq \overline{\mathcal{MH}^{m\times n}_{d,r,a}}$ immediately follows, which implies $\overline{\mathcal{MH}^{m\times n}_{d,r,a}} = \overline{\mathcal{B}^{m\times n}_{d,r,a}}$, which in turn implies the result in item (ii). Therefore, we focus on proving $\mathcal{B}^{m\times n}_{d,r,a} \subseteq \overline{\mathcal{MH}^{m\times n}_{d,r,a}}$.
If $L(\la) R(\la) \in \mathcal{B}^{m\times n}_{d,r,a}$, then
\begin{align} \label{eq.auxunderlinef}
  L(\la)^T & \in \CC[\la]^{r \times m}_{\underline{\mathbf{f}}}, \; \mbox{with $\underline{\mathbf{f}} = (\underbrace{d-d_R -1, \ldots, d-d_R - 1}_{t_R}  , \underbrace{d-d_R,\ldots, d-d_R}_{r-t_R})$,} \\ \label{eq.auxunderlineg}
  R(\la) & \in \CC[\la]^{r \times n}_{\underline{\mathbf{g}}}, \; \mbox{with $\underline{\mathbf{g}} = (\underbrace{d_R +1, \ldots, d_R + 1}_{t_R}  , \underbrace{d_R,\ldots, d_R}_{r-t_R})$.}
\end{align}
Therefore, Theorem \ref{thm.genericrobust} applied to $L(\la)$ and $R(\la)$ implies that there exist sequences of polynomial matrices $\{L_k(\la)\}_{k\in \NN}\subset \CC[\la]^{m\times r}$ and $\{R_k(\la)\}_{k\in \NN}\subset \CC[\la]^{r\times n}$, such that
\begin{itemize}
  \item[\rm (1)] $\lim_{k \rightarrow \infty} L_k(\la) = L(\la)$ and $\lim_{k \rightarrow \infty} R_k(\la) = R(\la)$,
  \item[\rm (2)] each polynomial matrix $L_k (\la)$ is a minimal basis, $\mathcal{N}_\ell (L_k)$ has minimal indices equal to
$\{\underbrace{\beta+1, \dots , \beta+1}_{t}, \underbrace{\beta, \dots , \beta}_{m-r-t}\}$, and $\deg ((L_k)_{*i}) = d - d_R - 1$ for $i=1,\ldots ,t_R$, and $\deg ((L_k)_{*i}) = d-d_R$ for $i=t_R + 1,\ldots ,r$,
\item[\rm (3)] each polynomial matrix $R_k (\la)$ is a minimal basis, $\mathcal{N}_r (R_k)$ has minimal indices equal to
$\{\underbrace{\alpha +1, \dots , \alpha +1}_{s}, \underbrace{\alpha, \dots , \alpha}_{n-r-s}\}$, and $\deg ((R_k)_{i*}) = d_R + 1$ for $i=1,\ldots ,t_R$, and $\deg ((R_k)_{i*}) = d_R$ for $i=t_R + 1,\ldots ,r$.
\end{itemize}
This means that $\{L_k(\la) R_k(\la)\}_{k\in \NN}\subset \mathcal{MH}^{m\times n}_{d,r,a}$ and that $\lim_{k \rightarrow \infty} L_k(\la) R_k (\la) = L(\la) R(\la)$. So, $L(\la) R(\la) \in \overline{\mathcal{MH}^{m\times n}_{d,r,a}}$ and $\mathcal{B}^{m\times n}_{d,r,a} \subseteq \overline{\mathcal{MH}^{m\times n}_{d,r,a}}$ is proved.

Items (iii) and (iv) follow from item (ii) and items (iii) and (iv) in Theorem \ref{thm.4descrition}.
\end{proof}

To compare the results we are obtaining for polynomial matrices with degree at most $d$, where $d \geq 1$, with those in \cite{DDL} for matrix pencils, that is, for $d=1$, we introduce some additional sets of polynomial matrices and prove for them a result similar to Theorem \ref{thm.4descrition}.

\begin{deff} \label{def.ineqhomogsubsets} Let $m,n,r$ and $d$ be integers such that $m,n \geq 2$, $d \geq 1$ and $0 < r < \min\{m,n\}$, and $a$ be an integer such that $0 \leq a \leq rd$. Let us define $d_R := \left\lfloor a /r \right\rfloor$, $t_R := a \; \mathrm{mod} \; r$, and the following subset of $\CC[\la]^{m\times n}_{d,r}$
$$
\mathcal{C}^{m\times n}_{d,r,a}  :=  \left\{ L(\la) R(\la) \, : \,
\begin{array}{l}
L(\la) \in \CC[\la]^{m\times r}, \; R(\la) \in \CC[\la]^{r\times n},\\
\deg (R_{i*}) \leq d_R + 1, \quad \mbox{for $i=1,\ldots ,t_R$}, \\
\deg (R_{i*}) \leq d_R, \quad \mbox{for $i=t_R + 1,\ldots ,r$}, \\
\deg (L_{*i}) \leq d - d_R -1, \quad \mbox{for $i=1,\ldots ,t_R$}, \\
\deg (L_{*i}) \leq d - d_R , \quad \mbox{for $i=t_R + 1,\ldots ,r$}
\end{array}
\right\}  \, .
$$
\end{deff}

\begin{theo} \label{thm.6descrition} Let $\mathcal{B}^{m\times n}_{d,r,a}$ and $\mathcal{C}^{m\times n}_{d,r,a}$ be the sets of polynomial matrices introduced in Definitions \ref{def.homogsubsets} and \ref{def.ineqhomogsubsets}, respectively. Then,
\begin{itemize}
\item[\rm (i)] $\mathcal{B}^{m\times n}_{d,r,a} \subseteq \mathcal{C}^{m\times n}_{d,r,a}$  for $a = 0,1,\ldots , rd$,

\item[\rm (ii)] $\overline{\mathcal{B}^{m\times n}_{d,r,a}} = \overline{\mathcal{C}^{m\times n}_{d,r,a}}$ for $a = 0,1,\ldots , rd$,

\item[\rm (iii)] $\displaystyle \CC[\la]^{m\times n}_{d,r} = \bigcup_{0 \leq a \leq rd} \overline{\mathcal{C}^{m\times n}_{d,r,a}} \, ,$ and

\item[\rm (iv)] for every $P(\la) \in \CC[\la]^{m\times n}_{d,r}$, there exists an integer $a$ such that $P(\la) \in \overline{\mathcal{C}^{m\times n}_{d,r,a}}$.
\end{itemize}
\end{theo}
\begin{proof}
Item (i) is obvious from the definitions of the involved sets.

Proof of item (ii). From item (i), we get that $\overline{\mathcal{B}^{m\times n}_{d,r,a}} \subseteq \overline{\mathcal{C}^{m\times n}_{d,r,a}}$. Next, we prove that $\mathcal{C}^{m\times n}_{d,r,a} \subseteq \overline{\mathcal{B}^{m\times n}_{d,r,a}}$. Let $L(\la) R(\la) \in \mathcal{C}^{m\times n}_{d,r,a}$, but $L(\la) R(\la) \notin \mathcal{B}^{m\times n}_{d,r,a}$. This means that the degrees of some rows of $R(\la)$ and/or of some columns of $L(\la)$ are strictly less than the corresponding quantities $d_R + 1$, $d_R$, $d-d_R - 1$, $d-d_R$ appearing in Definition \ref{def.ineqhomogsubsets}. Using any sequences of constant nonzero vectors $\{v_k\}_{k\in \NN} \subset \CC^{m\times 1}$ and/or $\{w_k\}_{k\in \NN} \subset \CC^{1\times n}$, such that $\lim_{k\rightarrow \infty} v_k = 0$ and $\lim_{k\rightarrow \infty} w_k = 0$, we sum to the rows of $R(\la)$ with degrees strictly less than $d_R + 1$ and/or $d_R$ polynomial vectors $\la^{d_R + 1} w_k$ and/or $\la^{d_R} w_k$, and sum to the columns of $L(\la)$ with degrees strictly less than $d-d_R - 1$ and/or $d-d_R$ polynomial vectors $\la^{d-d_R - 1} v_k$ and/or $\la^{d-d_R} v_k$. This allows us to construct a sequence $\{L_k(\la) R_k(\la)\}_{k\in \NN} \subset \mathcal{B}^{m\times n}_{d,r,a}$ such that $\lim_{k\rightarrow \infty} L_k(\la) R_k(\la) =
L (\la) R (\la)$. This proves $L (\la) R (\la) \in \overline{\mathcal{B}^{m\times n}_{d,r,a}}$ and $\mathcal{C}^{m\times n}_{d,r,a} \subseteq \overline{\mathcal{B}^{m\times n}_{d,r,a}}$, which implies $\overline{\mathcal{C}^{m\times n}_{d,r,a}} \subseteq \overline{\mathcal{B}^{m\times n}_{d,r,a}}$. This proves item (ii).

Items (iii) and (iv) follow from item (ii) and items (iii) and (iv) in Theorem \ref{thm.4descrition}.
\end{proof}

\begin{rem} \label{rem.compen1} {\rm (Comparisons with results for matrix pencils) For $d=1$, i.e., for matrix pencils, the sets $\mathcal{C}^{m\times n}_{1,r,a}$, for $a =0,1,\ldots , r$, in Definition \ref{def.ineqhomogsubsets} are exactly the sets $\mathcal{C}_a^r$ in \cite[Lemma 4]{DDL}. However, by using the Kronecker canonical form of pencils, Lemma 4 in \cite{DDL} proves that $\displaystyle \CC[\la]^{m\times n}_{1,r} = \bigcup_{0 \leq a \leq r} \mathcal{C}^{m\times n}_{1,r,a}$, which is a result stronger than Theorem \ref{thm.6descrition}-(iii) because it does not involve closures. This raises the question whether for $d\geq 2$ the closures can be removed in Theorem \ref{thm.6descrition}-(iii). Unfortunately, this is not possible as the next example shows.
}
\end{rem}

\begin{example} \label{ex.compen1} {\rm
Consider the polynomial matrix $P(\la) \in \CC[\la]^{3 \times 3}_{2,2}$ in \eqref{eq.extechnicalthm}. We are going to show that $P(\la) \notin \bigcup_{0 \leq a \leq 4} \mathcal{C}^{3\times 3}_{2,2,a}$. For this purpose, we follow an argument similar to that in Example \ref{ex.technicalthm}. Note first that the two factors $L(\la)$ and $R(\la)$ of $P(\la)$ in \eqref{eq.extechnicalthm} are minimal bases. Thus, the minimal indices of ${\mathcal Col} (P)$ are $2$ and $0$ and the  minimal indices of ${\mathcal Row} (P)$ are also $2$ and $0$. Moreover, since $\rank (P) = 2$, any factorization $P(\la) = \widetilde{L} (\la) \widetilde{R} (\la)$ with $\widetilde{L} (\la) \in \CC[\la]^{3\times 2}$ and $\widetilde{R} (\la) \in \CC[\la]^{2\times 3}$ must satisfy $\rank (\widetilde{L}) =\rank (\widetilde{R}) = 2$ and, so, the columns of $\widetilde{L} (\la)$ are a polynomial basis of ${\mathcal Col} (P)$ and the rows of $\widetilde{R} (\la)$ are a polynomial basis of ${\mathcal Row} (P)$. This means that the sum of the degrees of the columns of $\widetilde{L} (\la)$ must be larger than or equal to $2$ and that the sum of the degrees of the rows of $\widetilde{R} (\la)$ must be larger than or equal to $2$. Therefore, $P(\la) \notin \mathcal{C}^{3\times 3}_{2,2,0}$ and $P(\la) \notin \mathcal{C}^{3\times 3}_{2,2,1}$, because in both cases the sum of the degrees of the rows of $\widetilde{R} (\la)$ would be smaller than $2$, and also that $P(\la) \notin \mathcal{C}^{3\times 3}_{2,2,3}$ and $P(\la) \notin \mathcal{C}^{3\times 3}_{2,2,4}$, because in both cases the sum of the degrees of the columns of $\widetilde{L} (\la)$ would be smaller than $2$. Then, the only remaining option is $P(\la) \in \mathcal{C}^{3\times 3}_{2,2,2}$ but in this case $d_R = 1$ and $t_R = 0$, which implies that both rows of $\widetilde{R} (\la)$ must have degree exactly $1$, and that they will be a minimal basis of ${\mathcal Row} (P)$, which is impossible because the minimal indices of ${\mathcal Row} (P)$ are $2$ and $0$. {\blue \qed}}
\end{example}

\subsection{Relation between factorizations and generic complete eigenstructures in $\CC [\la]^{m\times n}_{d,r}$} \label{subsec.relation} A glance to the results in Theorems \ref{thm.mainthandrii2017}, \ref{thm.3description}, \ref{thm.4descrition}, \ref{thm.5descrition} and \ref{thm.6descrition} hints a relationship between the closures of the orbits $\orb (K_a)$ of polynomial matrices with generic eigenstructures and those of the sets defined before in Section \ref{sect.compact}. To establish this relationship, we characterize $\orb (K_a)$ as a set of factorized polynomial matrices in the next theorem.

\begin{theo} \label{thm.orbitsasfactors}
Let $m,n,r$ and $d$ be integers such that $m,n \geq 2$, $d \geq 1$ and $0 < r < \min\{m,n\}$, and $a$ be an integer such that $0 \leq a \leq rd$. Let us define $\alpha := \lfloor a / (n-r) \rfloor$, $s := a \mod (n-r)$, $\beta := \lfloor (rd-a)/(m-r) \rfloor$, and $t := (rd-a) \mod (m-r)$. Let $\orb (K_a)$ be the orbit of polynomial matrices in $\CC[\la]^{m\times n}_{d,r}$ appearing in Theorem \ref{thm.mainthandrii2017}. Then
$$
\orb (K_a)  =  \left\{ L(\la) R(\la) \, : \,
\begin{array}{l}
L(\la) \in \CC[\la]^{m\times r}, \; R(\la) \in \CC[\la]^{r\times n},\\
L(\la) \; \mbox{and} \; R(\la) \; \mbox{are minimal bases,}\\
\mathcal{N}_\ell (L) \; \mbox{has minimal indices} \;
\{\underbrace{\beta+1, \dots , \beta+1}_{t}, \underbrace{\beta, \dots , \beta}_{m-r-t}\} ,\\
\mathcal{N}_r (R) \; \mbox{has minimal indices} \;
\{\underbrace{\alpha+1, \dots , \alpha+1}_{s},\underbrace{\alpha, \dots , \alpha}_{n-r-s}\} ,\\
\deg (L_{*i}) + \deg (R_{i*}) = d , \quad \mbox{for $i=1,\ldots ,r$}
\end{array}
\right\}  \, .
$$
\end{theo}

\begin{proof} The result is an immediate corollary of Theorem \ref{thm.minfact} and the facts that $\mathcal{N}_\ell (L) = \mathcal{N}_\ell (L(\la) R(\la))$ and $\mathcal{N}_r (R) = \mathcal{N}_r (L(\la) R(\la))$ according to Lemma \ref{lemm:fact0}.
\end{proof}

With this result at hand, we get the next theorem.

\begin{theo} \label{thm.relationgeneigen} Let $\mathcal{A}^{m\times n}_{d,r,a}$, $\mathcal{B}^{m\times n}_{d,r,a}$, $\mathcal{M}^{m\times n}_{d,r,a}$,  $\mathcal{MH}^{m\times n}_{d,r,a}$ and $\mathcal{C}^{m\times n}_{d,r,a}$ be the sets of polynomial matrices introduced in Theorem \ref{thm.3description} and in Definitions \ref{def.homogsubsets}, \ref{def.minbassubsets} and \ref{def.ineqhomogsubsets}. Let $\orb (K_a)$ be the orbit of polynomial matrices in $\CC[\la]^{m\times n}_{d,r}$ appearing in Theorem \ref{thm.mainthandrii2017}. Then,
\begin{itemize}
\item[\rm (i)] $\orb (K_a) \subseteq \mathcal{A}^{m\times n}_{d,r,a}$  for $a = 0,1,\ldots , rd$,

\item[\rm (ii)] $\mathcal{MH}^{m\times n}_{d,r,a} \subseteq \orb (K_a)$ for $a = 0,1,\ldots , rd$,

\item[\rm (iii)] $\overline{\orb (K_a)} = \overline{\mathcal{MH}^{m\times n}_{d,r,a}}=\overline{\mathcal{M}^{m\times n}_{d,r,a}} = \overline{\mathcal{B}^{m\times n}_{d,r,a}} = \overline{\mathcal{C}^{m\times n}_{d,r,a}} = \overline{\mathcal{A}^{m\times n}_{d,r,a}}$ for $a = 0,1,\ldots , rd$.
\end{itemize}
\end{theo}

\begin{proof}
Proof of item (i). If $P(\la)\in \orb (K_a)$, then $P(\la) = L(\la) R(\la)$ with the factors $L(\la)$ and $R(\la)$ satisfying the properties described in Theorem \ref{thm.orbitsasfactors}. These properties imply that the degrees of the rows of $R(\la)$ are the minimal indices of ${\mathcal Row} (P)$ by Lemma \ref{lemm:fact0}. Combining this result with the fact that $\mathcal{N}_r (R) = \mathcal{N}_r (P)$, again by Lemma \ref{lemm:fact0}, and with Corollary \ref{cor.leftcolumnmin}, we get
\[
\sum_{i=1}^{r} \deg (R_{i*}) = s (\alpha + 1) + (n-r-s) \alpha = (n-r) \alpha + s = a.
\]
This implies that $P(\la)\in \mathcal{A}^{m\times n}_{d,r,a}$ and, so, item (i).

Item (ii) follows from the definitions of the involved sets.

Proof of item (iii). Item (i) implies $\overline{\orb (K_a)} \subseteq \overline{\mathcal{A}^{m\times n}_{d,r,a}}$. Combining this inclusion with Theorems \ref{thm.4descrition}-(ii), \ref{thm.5descrition}-(ii) and \ref{thm.6descrition}-(ii), we get
$$\overline{\orb (K_a)} \subseteq \overline{\mathcal{MH}^{m\times n}_{d,r,a}}=\overline{\mathcal{M}^{m\times n}_{d,r,a}} = \overline{\mathcal{B}^{m\times n}_{d,r,a}} = \overline{\mathcal{C}^{m\times n}_{d,r,a}} = \overline{\mathcal{A}^{m\times n}_{d,r,a}}.$$
On the other hand, item (ii) implies $\overline{\mathcal{MH}^{m\times n}_{d,r,a}} \subseteq \overline{\orb (K_a)}$, which combined with the equation above yields the result in item (iii).
\end{proof}

The inclusion relationships presented in Theorem \ref{thm.relationgeneigen}-(i) and (ii) between $\orb (K_a)$ and the other sets involved in this theorem are the only ones that hold in general. We illustrate this statement in the next example.

\begin{example} \label{ex.noinclusions} {\rm
Consider the following polynomial matrix
\begin{equation}\label{eq.noinclu1}
P(\la) =
\begin{bmatrix}
  0 & 0 & 1 & 1 \\
  0 & 0 & \la^2 & \la^2 \\
  1 & \la^2 & 2 \la^4 & \la^4 \\
  1 & \la^2 & \la^4 & 0
\end{bmatrix} =
\begin{bmatrix}
  1 & 0 \\
  \la^2 & 0 \\
  \la^4 & 1 \\
  0 & 1
\end{bmatrix}
\begin{bmatrix}
  0 & 0 & 1 & 1 \\
  1 & \la^2 & \la^4 & 0
\end{bmatrix} =: L(\la) R(\la).
\end{equation}
$P(\la)$ belongs to $\CC [\la]^{4\times 4}_{4,2}$. Moreover, the factors $L(\la)$ and $R(\la)$ are minimal bases by Theorem \ref{thm.minbasischar}. Consider also the following polynomial matrices
\begin{equation}\label{eq.noinclu2}
\widehat{L}(\la) = \left[ \begin{array}{rr}
  \la^2 & 0 \\
  -1 & \la^2 \\
  0 & -1 \\
  0 & 1
\end{array} \right] \quad \mbox{and} \quad
\widehat{R}(\la) =
\left[ \begin{array}{rrrr}
  \la^2 & -1 & 0 & 0 \\
  0 & \la^2 & -1 & 1
\end{array} \right] .
\end{equation}
It is easy to check that the columns of $\widehat{L}(\la)$ are a minimal basis of ${\mathcal N}_r (R) = {\mathcal N}_r (P)$ and that the rows of $\widehat{R}(\la)$ are a minimal basis of ${\mathcal N}_\ell (L) = {\mathcal N}_\ell (P)$. Therefore, $P(\la) \in \orb (K_4)$, by Theorem \ref{thm.orbitsasfactors}. However, $P(\la) \notin \mathcal{B}^{4\times 4}_{4,2,4}$, $P(\la) \notin \mathcal{M}^{4\times 4}_{4,2,4}$,  $P(\la) \notin \mathcal{MH}^{4\times 4}_{4,2,4}$ and $P(\la) \notin \mathcal{C}^{4\times 4}_{4,2,4}$. To see this, we need to check that no factorization of $P(\la)$ as $P(\la) = \widetilde{L}(\la) \widetilde{R} (\la)$, with $\widetilde{L}(\la) \in \CC[\la]^{4\times 2}$ and $\widetilde{R} (\la) \in \CC[\la]^{2\times 4}$, satisfies the conditions of the definitions of these sets. Note that in any of these factorizations $P(\la) = \widetilde{L}(\la) \widetilde{R} (\la)$ the rows of $\widetilde{R} (\la)$ are a polynomial basis of ${\mathcal Row} (P)$. Therefore, combining {\blue Theorem \ref{thm.mackeystrong}}
with the fact that the minimal indices of ${\mathcal Row} (P)$ are $0,4$, we obtain that $\deg (\widetilde{R}) \geq 4$. But, $d_R = \lfloor a/r \rfloor = \lfloor 4/2 \rfloor  =2$ and $t_R = 0$, which implies that any polynomial matrix in any of the sets $\mathcal{B}^{4\times 4}_{4,2,4}$, $\mathcal{M}^{4\times 4}_{4,2,4}$,  $\mathcal{MH}^{4\times 4}_{4,2,4}$ and $\mathcal{C}^{4\times 4}_{4,2,4}$ can be factorized as $L_S(\la) R_S (\la)$ with $L_S(\la)\in \CC[\la]^{4\times 2}$, $R_S(\la)\in \CC[\la]^{2\times 4}$ and $\deg (R_S) \leq 2$. Thus, $P(\la)$ does not belong to any of these sets.

Next, consider the polynomial matrix
\begin{equation}\label{eq.noinclu3}
Q(\la) =
 \left[ \begin{array}{cr}
  \la^2 & 0 \\
  -1 & \la^2 \\
  0 & -1 \\
  0 & 1
\end{array} \right]
\left[ \begin{array}{rrrr}
  \la^2 & -1 & 0 & 0 \\
  0 & \la^2 & -1 & 1
\end{array} \right]
= \left[ \begin{array}{rcrr}
  \la^4 & -\la^2 & 0 & 0\\
  -\la^2 & \la^4 +1 & -\la^2 & \la^2\\
  0 & -\la^2 & 1 &-1\\
  0 & \la^2 & -1& 1
\end{array} \right],
\end{equation}
which has been constructed as $Q(\la)= \widehat{L} (\la) \widehat{R} (\la)$ with the matrices in \eqref{eq.noinclu2}. Observe that $Q(\la) \in  \mathcal{M}^{4\times 4}_{4,2,4} \subseteq \mathcal{B}^{4\times 4}_{4,2,4} \subseteq \mathcal{C}^{4\times 4}_{4,2,4}$ and $Q(\la) \in \mathcal{A}^{4\times 4}_{4,2,4}$. However, $Q(\la) \notin  \orb (K_4)$ because the minimal indices of ${\mathcal N}_r (\widehat{R}) = {\mathcal N}_r (Q)$ are $0$ and $4$, since the columns of $L(\la)$ in \eqref{eq.noinclu1} are a minimal basis of ${\mathcal N}_r (\widehat{R})$. {\blue \qed}
}
\end{example}

\begin{rem} \label{rem.compen2} {\rm (Comparisons with results for matrix pencils) For $d = 1$, it was proved in \cite[Theorem 6]{DDL} that $\overline{\orb (K_a)} = \mathcal{C}^{m\times n}_{1,r,a}$, while Theorem \ref{thm.relationgeneigen} only proves the weaker result $\overline{\orb (K_a)} = \overline{\mathcal{C}^{m\times n}_{1,r,a}}$. For $d \geq 2$, the result $\overline{\orb (K_a)} = \overline{\mathcal{C}^{m\times n}_{d,r,a}}$  cannot be improved, since, in general, $\overline{\orb (K_a)} \ne \mathcal{C}^{m\times n}_{d,r,a}$. The polynomial matrix in \eqref{eq.noinclu1} illustrates this inequality.
}
\end{rem}

\section{Conclusions} \label{sec.conclusions}
We have established many results on rank factorizations and minimal rank factorizations of polynomial matrices, which, as far as we know, are completely new in the literature. In addition, the generic degree properties in the set $\CC[\la]^{m\times n}_{d,r}$ of complex $m\times n$ polynomial matrices of degree at most $d$ and rank at most $r$ of such factorizations have been carefully studied and several dense subsets of factorized polynomial matrices have been identified. Some of these subsets allow us to approximate any polynomial matrix in $\CC[\la]^{m\times n}_{d,r}$ as the limit of a sequence of factorized polynomial  matrices that can be easily and efficiently generated due to the particular degree properties of their factorizations, which have left factors with columns whose degrees differ at most by one and right factors with rows whose degrees differ at most by one. Apart from their fundamental nature in the theory of polynomial matrices, we hope that these results will have applications in the solution of different nearness problems involving polynomial matrices in $\CC[\la]^{m\times n}_{d,r}$. Possible lines of future research include exploring the development of structured rank factorizations and minimal rank factorizations of classes of structured polynomial matrices appearing in applications \cite{goodvibrations}, and verifying if some of the dense subsets of polynomial matrices in Section \ref{sect.compact} are also open in $\CC[\la]^{m\times n}_{d,r}$.

\bigskip \noindent
{\blue {\bf Acknowledgements.}
The authors thank two anonymous referees for many useful suggestions that have contributed to improve the presentation of this paper.}

\end{document}